\documentclass[12pt]{amsart}
\usepackage[utf8]{inputenc}
\usepackage{amsmath, amssymb, amsbsy, amsthm, amscd, amsfonts, latexsym, amstext, delarray, amsmath, graphicx, subcaption}
\usepackage{hyperref}
\usepackage{amssymb}
\usepackage{graphicx}
\usepackage{tikz}
\usepackage{dsfont}
\usepackage{mathtools}
\usepackage[nobysame, alphabetic]{amsrefs}
\usepackage[margin=1in]{geometry}
\usepackage[justification=centering]
{caption}
\usepackage{import}
\usepackage{xifthen}
\usepackage{pdfpages}
\usepackage{transparent}
\newtheorem{case}{Case}
\newtheorem{theorem}{Theorem}[section]

\newtheorem{lemma}[theorem]{Lemma}

\newtheorem{proposition}[theorem]{Proposition}
\newtheorem{conjecture}[theorem]{Conjecture}

\theoremstyle{definition}

\newtheorem{definition}[theorem]{Definition}

\title{Smooth Knots with Odd Quadratic Term of the Conway Polynomial have Inscribed Trefoils}
\author{Jonah Yoshida}
\date{August 2024}

\begin{document}
\begin{abstract}
    An \textit{inscribed knot} is formed by polygonally connecting points lying on a knot \(\gamma\) in parametric order, then closing the path by connecting the first and final points. The \textit{stick-knot number} of a knot type K is the minimum number of line segments needed to polygonally form some knot with the same homotopy type. The stick-knot number of a trefoil is 6. Cole Hugelmeyer studied the manifold \(M\) consisting of 6 points lying on a triangular prism and found that by intersecting a perturbation of \(M'\), twisting the top of the prism, with \(Q_\gamma\), the manifold of 6-tuples of points lying on \(\gamma\), any analytic knot with nontrivial quadratic term of its Conway polynomial has an inscribed trefoil. We show that by using a perturbation of the double-cover of the orientation class \([M \cap Q_\gamma]\) and analysis of planar configurations, an analogous result holds for a class of smooth knots with odd quadratic term. We also show that in the analytic case, there are both inscribed left and right-handed trefoils.
\end{abstract}
\maketitle
\section{Introduction}
\par Consider the generic knot embedded in \(\mathbb{R}^3\) and parametrized by \(\gamma: \mathbb{R} \rightarrow \mathbb{R}^3\), such that the function is analytic, \(\mathbb{Z}\)-periodic, with non-vanishing derivative (a regular curve), and of knot type \(K\). We can investigate whether there always exists a sequence of numbers \(0 \leq t_1 < t_2 < \ldots < t_6 < 1\) so that the polygonal path 
obtained by sequentially connecting the points \(\gamma(t_1), \gamma(t_2), \ldots, \gamma(t_6)\) by line segments is a trefoil. More generally, we can investigate whether there exist (potentially more) points so that connecting them in the same manner yields a different knot type. Since increasing the number of points in the sequentially-connected set of lines has a polygonal knot type approaching \(K\), there are only finitely many inscribed knot types. Another area of interest is determining the minimum non-zero crossing number of all inscribed knots for a nontrivial class of knots. In \cite{hugelmeyer}, Cole Hugelmeyer proved the following theorem.
\begin{theorem}
Let \(K\) be a knot type for which the quadratic term of the Conway polynomial is nontrivial, and let \(\gamma: \mathbb{R} \rightarrow \mathbb{R}^3\) be an analytic \(\mathbb{Z}\)-periodic function with nonvanishing derivative which parametrizes a knot of type \(K\) in space. Then, there exists a sequence of numbers \(0 \leq t_1 < t_2 < \ldots < t_6 < 1\) so that the polygonal path obtained by cyclically connecting the points \(\gamma(t_1), \gamma(t_2), \ldots, \gamma(t_6)\) by line segments is a trefoil knot.
\end{theorem}
In this paper, we prove the following two theorems.
\begin{theorem}
Let \(K\) be a knot type for which the quadratic term of the Conway polynomial is nontrivial, and let \(\gamma: \mathbb{R} \rightarrow \mathbb{R}^3\) be an analytic \(\mathbb{Z}\)-periodic function with nonvanishing derivative which parametrizes a knot of type \(K\) in space. Then, there exists a sequence of numbers \(0 \leq t_1 < t_2 < \ldots < t_6 < 1\) so that the polygonal path obtained by cyclically connecting the points \(\gamma(t_1), \gamma(t_2), \ldots, \gamma(t_6)\) by line segments is a left-handed trefoil knot. There also exists a sequence of numbers \(0 \leq t_1 < t_2 < \ldots < t_6 < 1\) so that the polygonal path obtained by cyclically connecting the points \(\gamma(t_1), \gamma(t_2), \ldots, \gamma(t_6)\) by line segments is a right-handed trefoil knot.
\end{theorem}
\begin{theorem}
Let \(K\) be a knot type for which the quadratic term of the Conway polynomial is odd, and let \(\gamma: \mathbb{R} \rightarrow \mathbb{R}^3\) be a smooth \(\mathbb{Z}\)-periodic function with nonvanishing derivative which parametrizes a knot of type \(K\) in space. Then, there exists a sequence of numbers \(0 \leq t_1 < t_2 < \ldots < t_6 < 1\) so that the polygonal path obtained by cyclically connecting the points \(\gamma(t_1), \gamma(t_2), \ldots, \gamma(t_6)\) by line segments is a trefoil knot.
\end{theorem}
\section{Left and Right-Handed Trefoils}
We begin by defining the manifold \(M\), the manifold of triangular prisms, which will appear in both following arguments.
\begin{definition}
    Let \(M\) be the moduli space for the following data:
\par 1) A point \(p \in \mathbb{R}\mathbb{P}^4 \backslash B^4\).
\par 2) An \textit{unordered} triple of distinct lines \(\{\ell_1,\ell_2,\ell_3\}\) that each pass through \(p\) and each intersect the 3-sphere \(\partial B^4\)
\par 3) A specified graph isomorphism \(G^c(\ell_1,\ell_2,\ell_3) \rightarrow \text{Cyc}_6\).
\end{definition}
To prove \textit{Theorem 1.2}, it suffices to show that analogies of Propositions 20, 21, and 22 of \cite{hugelmeyer} hold. The analogy we want for Proposition 20 is the same except for in the 1) clause,  we need a left or right-handed trefoil specified. 

\begin{proposition}
For some choice of \(M'\) as in Proposition 14 of \cite{hugelmeyer}, there exist real numbers \(t_1 < \ldots < t_6 < t_1 + 1\) so that \((\gamma(t_1), \ldots, \gamma(t_6)) \in Q_\gamma \cap M'\) and such that one of the three following possibilities holds.
\newline
\par 1) The 6-tuple \((\gamma(t_1), \ldots, \gamma(t_6))\) forms a right-handed trefoil.
\newline
\par 2) The 6-tuple \((\gamma(t_1), \ldots, \gamma(t_6))\) lies in some plane \(P\), is a stereographic trefoil, and at most three of the points in the 6-tuple are one-sided intersection points in \(P\).
\newline
\par 3) The 6-tuple \((\gamma(t_1), \ldots, \gamma(t_6))\) lies in some circle \(C\) of finite radius, and if \(P\) is the plane containing the circle \(C\), at most two of the points in the 6-tuple are one-sided in \(P\).
\end{proposition}
\textit{Proof.} Following the proof of Hugelmeyer's Proposition 20, it suffices to show that either \(Q_\gamma \cap M'\) consists of planar configurations or has a right-handed spherical trefoil. If we choose some orientation \(M'\) for the rotation of the top circle, assume for sake of contradiction that \(Q_\gamma \cap M'\) has its only spherical trefoils being left-handed. Then, choosing the opposite orientation for \(M'\) will yield right-handed trefoils since the parity of the handedness of trefoils is fixed with respect to a fixed point of inversion \(I\).
\(\boxed{}\)
\newline
\newline We proceed with the analogy of Proposition 21.
\begin{proposition}
If there exist real numbers \(t_1 < \ldots < t_6 < t_1 + 1\) so that the 6-tuple \((\gamma(t_1), \gamma(t_2), \ldots, \gamma(t_6))\) lies in some plane P, is a stereographic trefoil, and at most three of the points are one-sided intersection points in P, then there exist real numbers \(t_1' < \ldots < t_6' < t_1' + 1\) so that \((\gamma(t_1'), \gamma(t_2'), \ldots, \gamma(t_6'))\) gives a right-handed trefoil.
\end{proposition}
\textit{Proof. } Assume the case with three one-sided points, since all other cases can be reduced to this case. Observe that the set of trefoils is open in \(C_6(\mathbb{R}^3)\). Therefore, the set of stereographic trefoils is an open subset of the set of cospherical 6-tuples of points. We now find non-coplanar yet cospherical 6-tuples \((\gamma(t_1'), \ldots, \gamma(t_6'))\) arbitrarily close to \((\gamma(t_1), \ldots, \gamma(t_6))\). Choose an orientation for \(P\) and 
draw a circle \(C \in P\). Suppose first that the three one-sided points are not colinear. Then, let \(C\) be their circumcircle and \(\{S_i\}\) be a sequence of spheres of increasing radius with centers above \(P\) such that \(\forall i, C \in S_i\). By the openness of stereographic trefoils, there exists \(I\) such that \(\forall i>I\), the knots formed by connecting the points of \(S_i \cap \gamma\) are trefoils. Assuming for sake of contradiction that all such trefoils are left-handed trefoils, \(\forall i > I \hspace{1mm} \exists \delta\) such that each point of \(S_i \cap \gamma\) lies on a \(\delta\)-sphere centered on some point of \(P\). Let \(\{B_\delta\}\) be the collection of such spheres. Taking \(i\) sufficiently large, \(\delta \rightarrow 0\) such that \(\{B_\delta\} \cap \gamma \backslash S_i\) is a right-handed trefoil on \(\gamma\), since all of the crossings are flipped for this combination of points and no new ones are created with small enough \(\delta\), giving a homotopy-equivalent knot to a mirror image. This concludes the non-colinear case.
\par For the colinear case, consider each of the constructions in Lemmas 3.12, 3.13, and 3.14. The 3.14 construction can be directly applied, the 3.13 construction can be modified to take \(p_5'\) as some point arbitrarily-close to \(P\) and not on the line \(p_{3,4}\), and the 3.12 construction will have no bad configurations by assigning either of \(p_3\) or \(p_6\) some orthogonal height. \(\boxed{}\)
\begin{proposition}
Given six points of the knot \(\gamma\) on a circle, if at most two of them are one-sided relative to \(\gamma\), then we can form a right-handed trefoil on \(\gamma\). \end{proposition}
\textit{Proof.} 
Let \(\pi\) be the projection map onto the plane \(P\) of the circle given by an orientation of \(P\), and \(\sigma_\pi(x)\) the signed distance between \(x\) and \(\pi(x)\). Let \(\tilde{S} \subseteq S := \{\gamma(t_i)\}\) denote the two-sided points of \(S\). By \(\gamma\) analytic, \(\forall \epsilon\) \(\exists \delta\) such that \(\forall x \in S, y \in \gamma \hspace{1mm} \cap B_ \delta (x)\), \(|\sigma_\pi(y)| < \epsilon \leq \delta\) with \(\lim_{\epsilon} \delta(\epsilon) = 0\). Observe that for some fixed \(\delta\), \(\exists \epsilon_\text{min}=\inf_{x \in \tilde{S}, y \in \gamma \cap B_\delta(x)}|\sigma_{\pi}(y)|\) such that if \(x\) is a two-sided point, we can find \(y \in B_{\delta}(x)\) with \(\sigma_\pi(y)=h\), \(\forall h \in [-\epsilon_\text{min}, \epsilon_\text{min}]\). 
\par For sufficiently-small \(\delta\), if we can find a function \(f: S \rightarrow [-\epsilon_\text{min}, \epsilon_\text{min}]\), corresponding to \(\sigma_\pi(x)\) for six \(x\) each in a unique \(\delta\)-ball of elements of \(S\) with the correct crossing rules, then we will have found a right-handed trefoil. Since \(\gamma(t_i) \in Q_{\gamma} \cap M_0\), we have three possible combinations for where the \(\gamma(t_i)\) fall on a circle \(C\). Since all three are the same up to cyclic permutation, there is a unique case to find \(f\) for in which the three lines \(\ell_1, \ell_2, \ell_3\) all pass through \(C\) without intersecting each other. We describe different locations of the one-sided points for this case by considering subcases where \(U = S \backslash \tilde{S} \in \ker f\). It suffices to analyze each such case for \(|U| \leq 2\).
\par We can now formalize a possible set of valid right-handed trefoil crossing rules. If \(\ell_i\) is the line between \(\gamma(t_i)\) and \(\gamma(t_{i+1})\) \(\forall i<6\), and \(\ell_6\) is the line between \(\gamma(t_6)\) and \(\gamma(t_{1})\), then they are:
\begin{enumerate}
    \item \((1-\alpha_{11})f(\gamma(t_1))+\alpha_{11}f(\gamma(t_2)) < (1-\alpha_{41})f(\gamma(t_4))+\alpha_{41}f(\gamma(t_5))\)
    \item \((1-(\alpha_{11}+\alpha_{12}))f(\gamma(t_1))+(\alpha_{11}+\alpha_{12})f(\gamma(t_2)) < (1-\alpha_{34})f(\gamma(t_4))+\alpha_{34}f(\gamma(t_3))\)
    \item \((1-\alpha_{43})f(\gamma(t_5))+\alpha_{43}f(\gamma(t_4)) < (1-\alpha_{64})f(\gamma(t_1))+\alpha_{64}f(\gamma(t_6)) \)
    \item \((1-(\alpha_{64}+\alpha_{63}))f(\gamma(t_1))+(\alpha_{64}+\alpha_{63})f(\gamma(t_6)) < (1-(\alpha_{31}+\alpha_{32}))f(\gamma(t_3))+(\alpha_{31}+\alpha_{32})f(\gamma(t_4))\)
    \item \((1-\alpha_{31})f(\gamma(t_3))+\alpha_{31}f(\gamma(t_4)) < (1-\alpha_{51})f(\gamma(t_5))+\alpha_{51}f(\gamma(t_6))\)
    \item \((1-\alpha_{23})f(\gamma(t_3))+\alpha_{23}f(\gamma(t_2)) < (1-\alpha_{53})f(\gamma(t_6))+\alpha_{53}f(\gamma(t_5))\)
    \item \((1-\alpha_{21})f(\gamma(t_2))+\alpha_{21}f(\gamma(t_3)) < (1-\alpha_{61})f(\gamma(t_6))+\alpha_{61}f(\gamma(t_1))\), where \(\alpha_{ij}\) is the \(i^{\text{th}}\) line's \(j^{\text{th}}\) segment from \(\gamma(t_i)\) to \(\gamma(t_{i+1})\).
\newline
\newline
\end{enumerate}
\begin{case}
\(U \subset \{\gamma(t_1), \gamma(t_3), \gamma(t_4), \gamma(t_6)\}\)
\end{case}
Let \(\gamma(t_4)\) and \(\gamma(t_6)\) be the one-sided points. 
Now, letting \(f(\gamma(t_3)), f(\gamma(t_4)), f(\gamma(t_6)) = 0\), we can let \(f(\gamma(t_2))=+\epsilon, f(\gamma(t_1)) = -\epsilon',\) and \(f(\gamma(t_3)) = +\epsilon''\), such that \(\epsilon'' << \epsilon, \epsilon'\). Given these values, it is clear that everything on the line between \(\gamma(t_1)\) and \(\gamma(t_6)\) is below the plane \(P\) such that all crossings are resolved as desired. Then, with \(\epsilon'' << \epsilon\), the crossings involving the line between \(\gamma(t_5)\) and \(\gamma(t_6)\) are resolved. Finally, the line between \(\gamma(t_2)\) and \(\gamma(t_1)\) can pierce the triangle formed by \(\gamma(t_3), \gamma(t_4)\), and \(\gamma(t_5)\) with sufficient choice of \(-\epsilon', +\epsilon\). Therefore, we have found a right-handed trefoil for \(U \subset \{\gamma(t_3), \gamma(t_4), \gamma(t_6)\}\).
\newline
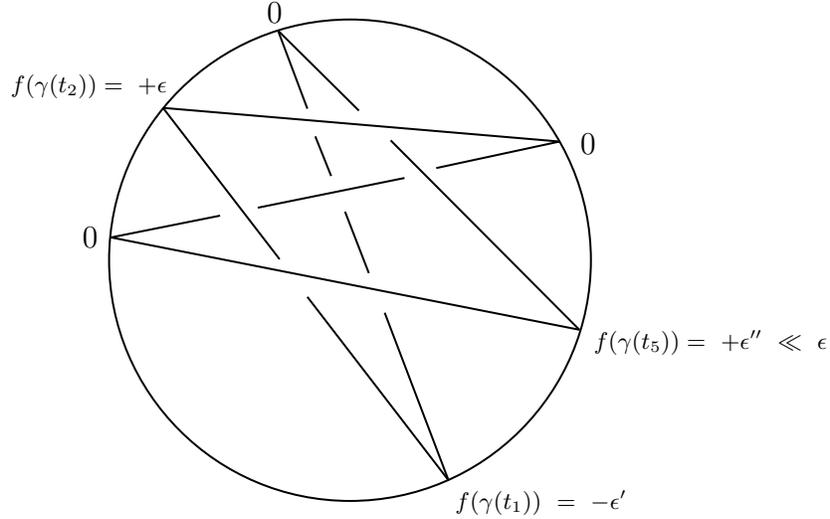
\begin{figure}
\begin{center}

\tikzset{every picture/.style={line width=0.75pt}} %set default line width to 0.75pt        

\begin{tikzpicture}[x=0.75pt,y=0.75pt,yscale=-1,xscale=1]

%Shape: Circle [id:dp1410965178983279] 
\draw   (123,150.5) .. controls (123,83.4) and (177.4,29) .. (244.5,29) .. controls (311.6,29) and (366,83.4) .. (366,150.5) .. controls (366,217.6) and (311.6,272) .. (244.5,272) .. controls (177.4,272) and (123,217.6) .. (123,150.5) -- cycle ;
%Straight Lines [id:da3725431767166054] 
\draw    (198,124) -- (272,109) ;
%Straight Lines [id:da4693334288826674] 
\draw    (150.2,73.6) -- (350.2,90.6) ;
%Straight Lines [id:da11503784760253666] 
\draw    (223,169) -- (294.2,261.6) ;
%Straight Lines [id:da04866695260393983] 
\draw    (227,88) -- (235,108) ;
%Straight Lines [id:da6138181910119653] 
\draw    (265,90) -- (360.2,185.6) ;
%Straight Lines [id:da721436217273137] 
\draw    (124,139) -- (360.2,185.6) ;
%Straight Lines [id:da17417925393698175] 
\draw    (208.2,34.6) -- (249,75) ;
%Straight Lines [id:da7562973365378978] 
\draw    (124,139) -- (179,128) ;
%Straight Lines [id:da31960302642683924] 
\draw    (288,104) -- (350.2,90.6) ;
%Straight Lines [id:da2564621511647762] 
\draw    (208.2,34.6) -- (223,74) ;
%Straight Lines [id:da532471790664923] 
\draw    (262,177) -- (294.2,261.6) ;
%Straight Lines [id:da1616798572978737] 
\draw    (242,126) -- (254,157) ;
%Straight Lines [id:da8613274628415968] 
\draw    (150.2,73.6) -- (209,150) ;

% Text Node
\draw (248,35) node   [align=left] {\begin{minipage}[lt]{68pt}\setlength\topsep{0pt}
0\\
\end{minipage}};
% Text Node
\draw (155,148) node   [align=left] {\begin{minipage}[lt]{68pt}\setlength\topsep{0pt}
0\\
\end{minipage}};
% Text Node
\draw (393,50) node   [align=left] {\begin{minipage}[lt]{68pt}\setlength\topsep{0pt}

\end{minipage}};
% Text Node
\draw (366,184.4) node [anchor=north west][inner sep=0.75pt]  [font=\scriptsize]  {$f( \gamma ( t_{5})) =\ +\epsilon ''\ \ll \ \epsilon $};
% Text Node
\draw (296.2,265) node [anchor=north west][inner sep=0.75pt]  [font=\scriptsize]  {$f( \gamma ( t_{1})) \ =\ -\epsilon '$};
% Text Node
\draw (406,101) node   [align=left] {\begin{minipage}[lt]{68pt}\setlength\topsep{0pt}
0\\
\end{minipage}};
% Text Node
\draw (72,55.4) node [anchor=north west][inner sep=0.75pt]  [font=\scriptsize]  {$f( \gamma ( t_{2})) =\ +\epsilon $};

\end{tikzpicture}

\end{center}
\caption{Case with \(U \subset \{\gamma(t_3), \gamma(t_4), \gamma(t_6)\}\).}
\end{figure}
\par Let \(\sigma_r = (1 6) (2 5) (3 4) \in S_6\). Then, replacing \(f\) with \(f_r = f \circ \sigma_r\) and choosing \(f\) as above, we may find necessary values of \(+\epsilon, -\epsilon', +\epsilon''\) such that we may find a right-handed trefoil for \(U \subset \{\gamma(t_1),\gamma(t_3),\gamma(t_4)\}\). Likewise, given \(\sigma_f = (1 3)(4 6) \in S_6\), replacing \(f\) with \(f_f = -f \circ \sigma_f\) yields a right handed-trefoil with appropriate values of \(+\epsilon, -\epsilon', +\epsilon''\) for \(U \subset \{\gamma(t_1),\gamma(t_4),\gamma(t_6)\}\). For the \(U \subset \{\gamma(t_1),\gamma(t_3),\gamma(t_6)\}\) case, take \(f_{fr} = -f \circ \sigma_f \circ \sigma_r\), which represents a combination of a rotation and a reflection. Since we were able to construct a right-handed trefoil given all 3-subsets of the 4-element set, \(|U| \leq 2\) means we are finished.
\begin{case}
\(U \subset \{\gamma(t_2), \gamma(t_4), \gamma(t_1)\}, \{\gamma(t_1), \gamma(t_4), \gamma(t_5)\}, \{\gamma(t_2), \gamma(t_6), \gamma(t_3)\}, \{\gamma(t_3), \gamma(t_5), \gamma(t_6)\}\)
\end{case}
\begin{figure}
\begin{center}

\begin{tikzpicture}[x=0.75pt,y=0.75pt,yscale=-1,xscale=1]

%Shape: Circle [id:dp1410965178983279] 
\draw   (123,150.5) .. controls (123,83.4) and (177.4,29) .. (244.5,29) .. controls (311.6,29) and (366,83.4) .. (366,150.5) .. controls (366,217.6) and (311.6,272) .. (244.5,272) .. controls (177.4,272) and (123,217.6) .. (123,150.5) -- cycle ;
%Straight Lines [id:da3725431767166054] 
\draw    (198,124) -- (272,109) ;
%Straight Lines [id:da4693334288826674] 
\draw    (150.2,73.6) -- (350.2,90.6) ;
%Straight Lines [id:da11503784760253666] 
\draw    (223,169) -- (294.2,261.6) ;
%Straight Lines [id:da04866695260393983] 
\draw    (227,88) -- (235,108) ;
%Straight Lines [id:da6138181910119653] 
\draw    (265,90) -- (360.2,185.6) ;
%Straight Lines [id:da721436217273137] 
\draw    (124,139) -- (360.2,185.6) ;
%Straight Lines [id:da17417925393698175] 
\draw    (208.2,34.6) -- (249,75) ;
%Straight Lines [id:da7562973365378978] 
\draw    (124,139) -- (179,128) ;
%Straight Lines [id:da31960302642683924] 
\draw    (288,104) -- (350.2,90.6) ;
%Straight Lines [id:da2564621511647762] 
\draw    (208.2,34.6) -- (223,74) ;
%Straight Lines [id:da532471790664923] 
\draw    (262,177) -- (294.2,261.6) ;
%Straight Lines [id:da1616798572978737] 
\draw    (242,126) -- (254,157) ;
%Straight Lines [id:da8613274628415968] 
\draw    (150.2,73.6) -- (209,150) ;

% Text Node
\draw (157,152) node   [align=left] {\begin{minipage}[lt]{68pt}\setlength\topsep{0pt}
0\\
\end{minipage}};
% Text Node
\draw (188,77) node   [align=left] {\begin{minipage}[lt]{68pt}\setlength\topsep{0pt}
0\\
\end{minipage}};
% Text Node
\draw (368,183.4) node [anchor=north west][inner sep=0.75pt]  [font=\scriptsize]  {$f( \gamma ( t_{5})) =\ +\epsilon '$};
% Text Node
\draw (344.2,286.6) node   [align=left] {\begin{minipage}[lt]{68pt}\setlength\topsep{0pt}
0\\
\end{minipage}};
% Text Node
\draw (358,85.4) node [anchor=north west][inner sep=0.75pt]  [font=\scriptsize]  {$ \begin{array}{l}
f( \gamma ( t_{3})) =-\epsilon \\
\end{array}$};
% Text Node
\draw (162,13.4) node [anchor=north west][inner sep=0.75pt]  [font=\scriptsize]  {$f( \gamma ( t_{6})) =\ -\epsilon ''$};
% Text Node
\draw (378,102.4) node [anchor=north west][inner sep=0.75pt]  [font=\scriptsize]  {$ \begin{array}{l}
\epsilon \ \ll \ \epsilon ''\\
\end{array}$};

\end{tikzpicture}

\end{center}
\caption{Case with \(U \subset \{\gamma(t_1), \gamma(t_2), \gamma(t_4)\}\).}
\end{figure}
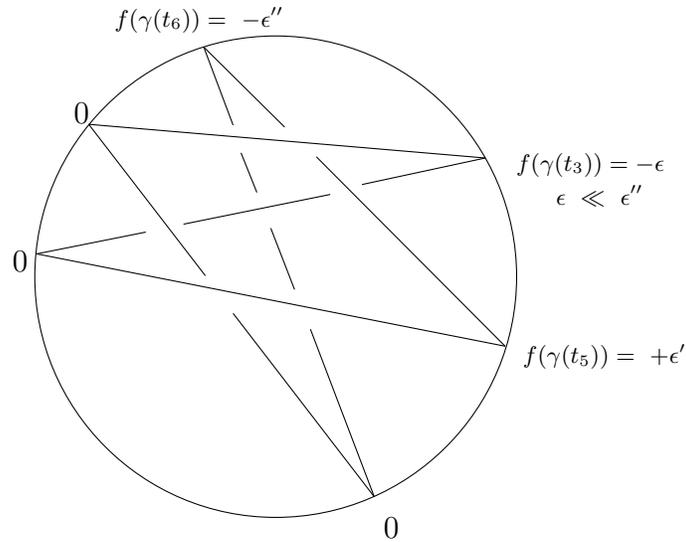
\par Let \(\gamma(t_1), \gamma(t_2),\) and \(\gamma(t_4)\) be the one-sided points. Then, given \(f(\gamma(t_1)) = f(\gamma(t_2)) = f(\gamma(t_4)) = 0\), let \(f(\gamma(t_6)) = -\epsilon''\), \(f(\gamma(t_3))=-\epsilon\), and \(f(\gamma(t_5)) = +\epsilon'\) such that \(\epsilon << \epsilon''\). Observe that all crossings on the line between \(\gamma(t_6)\) and \(\gamma(t_1)\) are resolved as \(\epsilon << \epsilon''\) and the strand from \(\gamma(t_4)\) to \(\gamma(t_5)\) is above \(P\). On the line between \(\gamma(t_1)\) and \(\gamma(t_2)\), crossings are resolved appropriately given the signs of \(f(\gamma(t_3)), f(\gamma(t_5))\). Lastly, since we can choose \(\epsilon\) arbitrarily close to \(0\), we can pick \(\epsilon', \epsilon''\) so that the line between \(\gamma(t_5)\) and \(\gamma(t_6)\) pierces the triangle formed between \(\gamma(t_2), \gamma(t_3),\) and \(\gamma(t_4)\). Thus, the right-handed trefoil. Applying the same \(\sigma_{r,f}\) operators gives the other 3-tuple of \(\gamma(t_i)\) cases.
\begin{case}
\(U \subset \{\gamma(t_2), \gamma(t_5)\}\)
\end{case}
Let \(f(\gamma(t_2))=f(\gamma(t_5))=0\), \(f(\gamma(t_3))=f(\gamma(t_4))=-\epsilon\), and \(f(\gamma(t_1))=f(\gamma(t_6))=-\epsilon'\) for \(\epsilon < \epsilon'\). Then, observe that the crossings along the segment from \(\gamma(t_1)\) to \(\gamma(t_2)\) may be resolved by picking sufficient \(\epsilon'\), since the segment from \(\gamma(t_3)\) to \(\gamma(t_4)\) is lower with respect to \(P\) than the segment from \(\gamma(t_4)\) to \(\gamma(t_5)\) along the original segment. Likewise for all crossings along the segment from \(\gamma(t_5)\) to \(\gamma(t_6)\). Finally, the crossings along the segment from \(\gamma(t_1)\) to \(\gamma(t_6)\) are resolved by \(-\epsilon'\) being the lowest orthogonal heights reached on the polygonal knot. Thus, the right-handed trefoil.
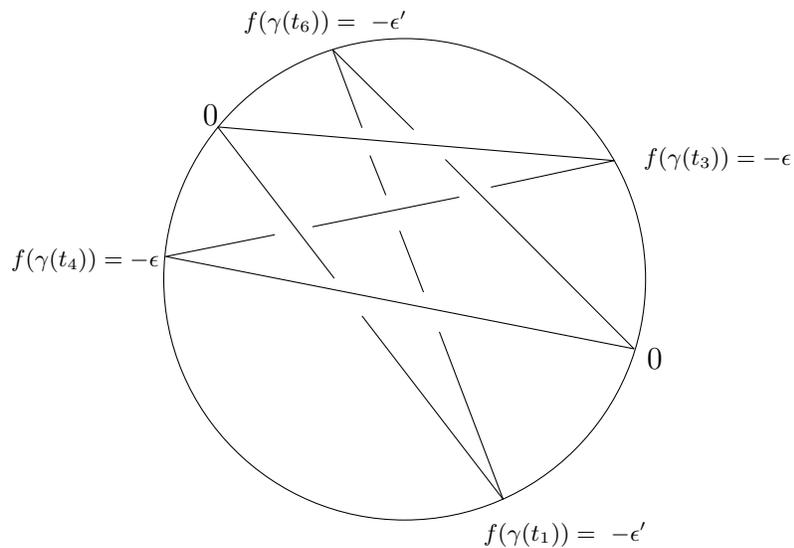
\begin{figure}
\begin{center}
\begin{tikzpicture}[x=0.75pt,y=0.75pt,yscale=-1,xscale=1]

%Shape: Circle [id:dp1410965178983279] 
\draw   (123,150.5) .. controls (123,83.4) and (177.4,29) .. (244.5,29) .. controls (311.6,29) and (366,83.4) .. (366,150.5) .. controls (366,217.6) and (311.6,272) .. (244.5,272) .. controls (177.4,272) and (123,217.6) .. (123,150.5) -- cycle ;
%Straight Lines [id:da3725431767166054] 
\draw    (198,124) -- (272,109) ;
%Straight Lines [id:da4693334288826674] 
\draw    (150.2,73.6) -- (350.2,90.6) ;
%Straight Lines [id:da11503784760253666] 
\draw    (223,169) -- (294.2,261.6) ;
%Straight Lines [id:da04866695260393983] 
\draw    (227,88) -- (235,108) ;
%Straight Lines [id:da6138181910119653] 
\draw    (265,90) -- (360.2,185.6) ;
%Straight Lines [id:da721436217273137] 
\draw    (124,139) -- (360.2,185.6) ;
%Straight Lines [id:da17417925393698175] 
\draw    (208.2,34.6) -- (249,75) ;
%Straight Lines [id:da7562973365378978] 
\draw    (124,139) -- (179,128) ;
%Straight Lines [id:da31960302642683924] 
\draw    (288,104) -- (350.2,90.6) ;
%Straight Lines [id:da2564621511647762] 
\draw    (208.2,34.6) -- (223,74) ;
%Straight Lines [id:da532471790664923] 
\draw    (262,177) -- (294.2,261.6) ;
%Straight Lines [id:da1616798572978737] 
\draw    (242,126) -- (254,157) ;
%Straight Lines [id:da8613274628415968] 
\draw    (150.2,73.6) -- (209,150) ;

% Text Node
\draw (412,200) node   [align=left] {\begin{minipage}[lt]{68pt}\setlength\topsep{0pt}
0\\
\end{minipage}};
% Text Node
\draw (188,77) node   [align=left] {\begin{minipage}[lt]{68pt}\setlength\topsep{0pt}
0\\
\end{minipage}};
% Text Node
\draw (162,13.4) node [anchor=north west][inner sep=0.75pt]  [font=\scriptsize]  {$f( \gamma ( t_{6})) =\ -\epsilon '$};
% Text Node
\draw (38,132.4) node [anchor=north west][inner sep=0.75pt]  [font=\scriptsize]  {$ \begin{array}{l}
f( \gamma ( t_{4})) =-\epsilon \\
\\
\end{array}$};
% Text Node
\draw (357,81.4) node [anchor=north west][inner sep=0.75pt]  [font=\scriptsize]  {$ \begin{array}{l}
f( \gamma ( t_{3})) =-\epsilon \\
\\
\end{array}$};
% Text Node
\draw (283,271.4) node [anchor=north west][inner sep=0.75pt]  [font=\scriptsize]  {$f( \gamma ( t_{1})) =\ -\epsilon '$};

\end{tikzpicture}
\end{center}
\caption{Case with \(U \subset \{\gamma(t_2), \gamma(t_5)\}\).}
\end{figure}

Therefore, we have considered all possible cases for two one-sided points, and have proven the proposition. \(\square\)
\newline
\par We now have a proof of Theorem 1.2.

\textit{Proof.} Combine propositions 2.2, 2.3, and 2.4. Then, repeat all arguments for left-handed trefoils. \(\square\)
\section{Trefoils on Smooth \(\gamma\)}
\par To prove the analogy of Theorem 1 
\cite{hugelmeyer} for a class of smooth non-trivial knots, we must instead consider the \(M\) configuration space with respect to only their cyclic ordering. 
\begin{definition}
    Let \(A\) be the group action on \(M\) which sends \(p_i\) to \(p_{i+1 \mod 6}\). Then, \(A \circlearrowright M\) induces equivalence classes on \(M\): \(m_1 \sim^{A} m_2\) if the elements have identical \(p, \ell_1, \ell_2, \ell_3\) and, if \(i_1, i_2\) are the graph isomorphisms for \(m_1, m_2\) respectively, \(c^j i_1 = i_2\) for some \(j \in \mathbb{Z}\) and \(c: \text{Cyc}_6 \rightarrow \text{Cyc}_6\) with \(c(i) = i+1\), corresponding to a shift of the ordering.
\end{definition}

\begin{lemma}
    Let \(\tilde{M} = M / A\). Then, \(H_1(\tilde{M}) \simeq \mathbb{Z}/6\mathbb{Z}\).
\end{lemma}
\textit{Proof.} Note that \(M\) is the universal covering space for \(\tilde{M}\). Thus, the abelianization of \(\pi_1(\tilde{M})\simeq \text{Deck}(M/\tilde{M}) \simeq \mathbb{Z}/6\mathbb{Z}\) is isomorphic to \(H_1(\tilde{M})\). Cyclic groups are abelian, so it suffices to show that \(\langle c^i \rangle\) generate the unique deck transformations of \(\tilde{M}\). Suppose there were a different \(g \in \text{Aut}(\text{Cyc}_6)\)
that generates a deck transformation. As \(\text{Aut}(\text{Cyc}_6) \simeq \mathbb{Z}/2\mathbb{Z}\), the only nontrivial automorphism \(g\) has \(g(0)=0, g(1)=5, g(2)=4, g(3)=3, g(4)=2,\) and \(g(5)=1\), corresponding to fixing \(p_1\) and \(p_4\) but swapping pairs \((p_2,p_6), (p_3,p_5)\). This reverses the cyclic ordering of the points and therefore cannot induce a deck transformation. \(\square\)
\newline
\par Observe that by the Gysin sequence and Poincaré duality, we have the following homomorphism induced by the pullback map from the normal bundle \(\tilde{M}_n\) to \(\tilde{M}\):
$$ \pi^{*} : H_1(\tilde{M}) \rightarrow H_5(\tilde{M}_n). $$
\par As a quick digression, it is worth noting that there may exist a proof of Theorem 1.3 by checking whether the homology class of \(Q_\gamma\) exists in the intersection of \(\pi^*(g)\) with the set of trefoils, for \(g\) a generator of \(H_1(\tilde{M})\) that forms a tubular neighborhood under \(\pi^*\).
\begin{proposition}
    If \(\gamma\) is a trefoil knot, then \(\kappa(\gamma) = [\tilde{M} \cap Q_\gamma] \in H_1(Q_{\gamma})\) is nontrivial.
\end{proposition}
\textit{Proof.} By Proposition 8 of \cite{hugelmeyer}, it suffices to show that for the parametric trefoil \(\gamma(t) = \frac{1}{\sqrt{2}}(\cos{4\pi t},\sin{4\pi t},\cos{6\pi t}, \sin{6\pi t})\), \(x \in Q_{\gamma}^p \implies x = ((\gamma(0), \gamma(\frac{1}{6}), \ldots, \gamma(\frac{5}{6})) \) and \(|Q_\gamma^{p} \cap \tilde{M}| = 1\). Let \(\pi_{12}\) be the projection of \(\gamma\) into its first two coordinates and \(\pi_{34}\) be the projection of \(\gamma\) into its last two coordinates. Since \(\pi_{12}(\gamma), \pi_{34}(\gamma) \subseteq S^1\), if \(x = (x_1, x_2, x_3, x_4, x_5, x_6) \in Q_{\gamma}^p \cap \tilde{M}\), then either \(\pi_{12}(x)\) is a set of points that lie on three distinct lines determined by their graph isomorphism and \(S^1\) or on three distinct points such that pairs \((x_1,x_4), (x_2,x_5), (x_3,x_6)\) are mapped to the same point under \(\pi_{12}\) and \(S^1\). Assume the former. Then, there exists some \((t_1, t_2, t_3, t_4, t_5, t_6)\) such that \(\forall t_i \in [0,1]\) and \(x_i\) is \((\text{Re}(e^{4 \pi t_i i}),\text{Im}(e^{4 \pi t_i i}))\) under \(\pi_{12}\). Thus, there exists a way of picking a point on \(S^1\) and traversing \(S^1\) counterclockwise so that we find the \(x_i\) in order and fewer than \(4 \pi\) radians elapse between encountering \(x_1\) and \(x_6\). However, by exhaustively checking all possible combinations, this cannot be done.
\begin{figure}
\begin{center}
\tikzset{every picture/.style={line width=0.75pt}} %set default line width to 0.75pt        

\begin{tikzpicture}[x=0.75pt,y=0.75pt,yscale=-1,xscale=1]

%Shape: Circle [id:dp3242021172293328] 
\draw   (79,124.43) .. controls (79,94.37) and (103.37,70) .. (133.43,70) .. controls (163.49,70) and (187.86,94.37) .. (187.86,124.43) .. controls (187.86,154.49) and (163.49,178.86) .. (133.43,178.86) .. controls (103.37,178.86) and (79,154.49) .. (79,124.43) -- cycle ;
%Straight Lines [id:da8526610179232306] 
\draw    (168.86,82.57) -- (173.86,160.57) ;
%Straight Lines [id:da0476787266828933] 
\draw    (148.86,71.57) -- (109.86,172.57) ;
%Straight Lines [id:da3018541403794268] 
\draw    (79,124.43) -- (119.86,71.57) ;
%Shape: Free Drawing [id:dp6102274976771969] 
\draw  [line width=3] [line join = round][line cap = round] (168.86,82.57) .. controls (168.86,82.57) and (168.86,82.57) .. (168.86,82.57) ;
%Shape: Free Drawing [id:dp842733138776121] 
\draw  [line width=3] [line join = round][line cap = round] (147.86,72.57) .. controls (147.86,72.57) and (147.86,72.57) .. (147.86,72.57) ;
%Shape: Free Drawing [id:dp09105878609924978] 
\draw  [line width=3] [line join = round][line cap = round] (173.86,160.57) .. controls (173.86,160.57) and (173.86,160.57) .. (173.86,160.57) ;
%Shape: Free Drawing [id:dp935335754306833] 
\draw  [line width=3] [line join = round][line cap = round] (109.86,172.57) .. controls (109.86,172.57) and (109.86,172.57) .. (109.86,172.57) ;
%Shape: Free Drawing [id:dp7958167569918486] 
\draw  [line width=3] [line join = round][line cap = round] (118.86,71.57) .. controls (118.86,71.57) and (118.86,71.57) .. (118.86,71.57) ;
%Shape: Free Drawing [id:dp09364117666154548] 
\draw  [line width=3] [line join = round][line cap = round] (78.86,124.57) .. controls (78.86,124.57) and (78.86,124.57) .. (78.86,124.57) ;
%Curve Lines [id:da04430542872708032] 
\draw    (195.86,74.43) .. controls (195.86,64.68) and (195.86,57.78) .. (174.54,58.37) ;
\draw [shift={(172.86,58.43)}, rotate = 357.51] [color={rgb, 255:red, 0; green, 0; blue, 0 }  ][line width=0.75]    (10.93,-3.29) .. controls (6.95,-1.4) and (3.31,-0.3) .. (0,0) .. controls (3.31,0.3) and (6.95,1.4) .. (10.93,3.29)   ;

% Text Node
\draw (173,74) node [anchor=north west][inner sep=0.75pt]    {$x_{1}$};
% Text Node
\draw (101,177.4) node [anchor=north west][inner sep=0.75pt]    {$x_{2}$};
% Text Node
\draw (109,55) node [anchor=north west][inner sep=0.75pt]    {$x_{3}$};
% Text Node
\draw (170,166.4) node [anchor=north west][inner sep=0.75pt]    {$x_{4}$};
% Text Node
\draw (145,57) node [anchor=north west][inner sep=0.75pt]    {$x_{5}$};
% Text Node
\draw (59,122.4) node [anchor=north west][inner sep=0.75pt]    {$x_{6}$};

\end{tikzpicture}
\caption{There is no counterclockwise path that passes through the points in the right order and elapses less than \(4\pi\) radians.}
\end{center}

\end{figure}
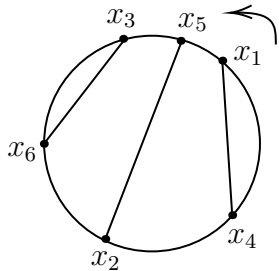
\par Hence, \(\pi_{12}(x)\) is a set of three unique points and \(|t_1-t_4|=|t_2-t_5|=|t_3-t_6|=2\pi\) such that the map \(T: t_i \mapsto \frac{3}{2} t_i\) 
induces a map \(T^*:\pi_{12}(\gamma) \rightarrow \pi_{34}(\gamma)\) where pairs \((x_1,x_4), (x_2, x_5), (x_3,x_6)\) are taken to antipodal points on \(S^1\). The lines \(\ell_1, \ell_2, \ell_3\) between \((x_1, x_4), (x_2,x_5),\) and \((x_3,x_6)\), respectively, have the forms \((a_1, b_1, c_1s, d_1s), (a_2, b_2, c_2s, d_2s), (a_3, b_3, c_3s, d_3s)\), for some parameter \(s\) and \(\forall \pi_{12}(x_i)=(a_i,b_i)\) and \(T(\pi_{12}(x_i))=(c_i,d_i) \in S^1\). Note that if \(\ell_1 \cap \ell_2 \cap \ell_3 \neq \emptyset\), then the angle formed between \(T^*(\pi_{12}(x_i))\) and \(y=0\) must be identical \(\forall i\), such that \(T(t_i)\) must be equally-spaced along with \(t_i\). Therefore, \(x=(\gamma(0), \gamma(\frac{1}{6}), \ldots, \gamma(\frac{5}{6}))\). \(\boxed{}\)
\newline
\par We now have tools to explore the smooth topological properties of \(\gamma\).
\begin{definition}
Let \(D\) be the double cover of \(\tilde{M}\) with respect to the orientation covering \(m \in \tilde{M} \mapsto (m,o)\), where \(o\) is an orientation on the circle containing the top triangle.
\end{definition}
\begin{lemma}
If \(\gamma\) has a nontrivial quadratic term in its Conway polynomial, there exists a stereographic trefoil on \(\gamma\) within \([D' \cap Q_\gamma]\), where \(D'\) is a perturbation of \(D\).
\end{lemma}

\textit{Proof.} By Proposition 3.3, \([\tilde{M} \hspace{0.5mm} \cap Q_{\gamma}]\) is a rank 2 Vassiliev invariant.  Observe that 
\newline
\(|H_1(\tilde{M}):H_1(D)| = 1,2\) and an element of \(H_1(\tilde{M})\) consists of a rotation and inversion of the triangle, reversing orientation, such that applying Lemma 3.2 gives \(H_1(D) \simeq \mathbb{Z}/3\mathbb{Z}\). Since \([D \hspace{0.5mm} \cap Q_{\gamma}] \in H_1(D)\) is a non-trivial orientable 1-manifold, we may perturb \(D\) to \(D'\) by rotating the oriented circle given by \([D \hspace{0.5mm} \cap Q_{\gamma}]\) along the oriented direction such that \([D' \hspace{0.5mm} \cap Q_{\gamma}] \subseteq S\), the set of stereographic trefoils. Thus, there exists a stereographic trefoil on \(\gamma\) which lies in the \([D' \hspace{0.5mm} \cap Q_{\gamma}]\) class. \(\square{}\)
\newline
\begin{lemma}
If the quadratic term of the Conway polynomial of \(\gamma\) is odd, then either \(\gamma\) has an inscribed trefoil or the stereographic projection of \([D' \cap Q_\gamma] \subset S\) lies within at least two planes.
\end{lemma}
\textit{Proof.} By Lemma 3.5, there exists a stereographic trefoil on \(\gamma\) in the \([D' \cap Q_\gamma]\) class. Let \(S_I\) denote the stereographic projection by inversion point \(I\). Assume for sake of contradiction that \(S_I([D' \cap Q_{\gamma}])\) lies entirely within a plane \(P\). Otherwise, we have either \(S_I([D' \cap Q_\gamma])\) consisting of at least 2 planes or a trefoil on \(\gamma\) by Proposition 13 of \cite{hugelmeyer}, either of which satisfies the conclusion. There is thus a continuous map \(p\) taking all points of \([D' \cap Q_\gamma]\) to planar configurations: \(p: [D' \cap Q_\gamma] \rightarrow \mathbb{R}\mathbb{P}^2\) with the induced homomorphism between first homologies given by the Hurewicz homomorphism \(p^*: H_1(D') \simeq \mathbb{Z} / 6 \mathbb{Z} \rightarrow H_1(\mathbb{R}\mathbb{P}^2) \simeq \mathbb{Z} / 2 \mathbb{Z}\). As the quadratic term of the Conway polynomial of \(\gamma\) is odd, any nontrivial homomorphism \(p^*\) will take \(g=[D' \cap Q_\gamma] \in H_1(D')\) to a nontrivial element of \(H_1(\mathbb{R}\mathbb{P}^2)\), a contradiction of \(S_I([D' \cap Q_\gamma])\) lying entirely within \(P\). So, \(S_I([D' \cap Q_{\gamma}])\) lies within at least \(2\) distinct planes \(P_1, P_2\). \(\square\)
\newline
\par If we can construct a trefoil given the latter condition of two or more planes containing \(S_I([D' \cap Q_\gamma])\), then we will have a proof of Theorem 1.3.
\par Assuming \(S_I([D' \cap Q_\gamma])\) consists entirely of planar configurations, or else there is a trefoil on \(\gamma\), we may categorize such configurations as shown in Figure 5. The first configuration type represents an identification of opposite points of \(\mathbb{R}\mathbb{P}^4\), the second a transition where the point of inversion \(I\) lies on the circumcircle of points, the third when the rotation is halfway complete, the fourth when complete, and the fifth after a planar flip induced by moving the configuration in \(S^3\). We will now prove a few lemmas that will be helpful for reducing cases for each of these configurations types when three points are colinear.

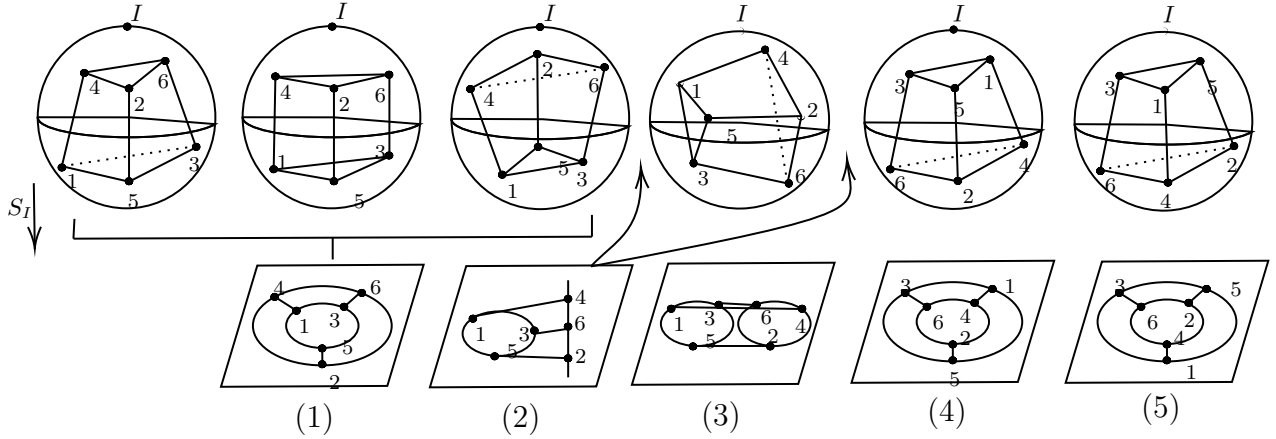
\begin{figure}
\begin{center}

\tikzset{every picture/.style={line width=0.75pt}} %set default line width to 0.75pt        

\begin{tikzpicture}[x=0.75pt,y=0.75pt,yscale=-1,xscale=1]
%uncomment if require: \path (0,299); %set diagram left start at 0, and has height of 299

%Shape: Ellipse [id:dp028681051683502323] 
\draw   (21.98,74.71) .. controls (21.98,49.32) and (42.03,28.74) .. (66.78,28.74) .. controls (91.52,28.74) and (111.58,49.32) .. (111.58,74.71) .. controls (111.58,100.1) and (91.52,120.68) .. (66.78,120.68) .. controls (42.03,120.68) and (21.98,100.1) .. (21.98,74.71) -- cycle ;
%Shape: Arc [id:dp9312617551658195] 
\draw  [draw opacity=0] (111.51,78.19) .. controls (104.96,82.03) and (88.05,84.79) .. (68.21,84.83) .. controls (42.68,84.88) and (21.98,80.42) .. (21.97,74.86) .. controls (21.97,74.81) and (21.97,74.76) .. (21.97,74.71) -- (68.19,74.76) -- cycle ; \draw   (111.51,78.19) .. controls (104.96,82.03) and (88.05,84.79) .. (68.21,84.83) .. controls (42.68,84.88) and (21.98,80.42) .. (21.97,74.86) .. controls (21.97,74.81) and (21.97,74.76) .. (21.97,74.71) ;  
%Shape: Free Drawing [id:dp07230902753389179] 
\draw  [line width=3] [line join = round][line cap = round] (45.29,52.38) .. controls (45.71,52.38) and (45.71,52.38) .. (45.29,52.38) ;
%Shape: Free Drawing [id:dp34672837808122914] 
\draw  [line width=3] [line join = round][line cap = round] (67,28.96) .. controls (67.42,28.96) and (67.42,28.96) .. (67,28.96) ;
%Shape: Free Drawing [id:dp17983411634163882] 
\draw  [line width=3] [line join = round][line cap = round] (33.99,100.08) .. controls (34.42,100.08) and (34.42,100.08) .. (33.99,100.08) ;
%Shape: Free Drawing [id:dp23731192697501768] 
\draw  [line width=3] [line join = round][line cap = round] (86.11,46.31) .. controls (86.53,46.31) and (86.53,46.31) .. (86.11,46.31) ;
%Shape: Free Drawing [id:dp259902022710234] 
\draw  [line width=3] [line join = round][line cap = round] (67.87,60.18) .. controls (68.29,60.18) and (68.29,60.18) .. (67.87,60.18) ;
%Straight Lines [id:da21158275797827364] 
\draw    (34.15,99.43) -- (45.44,52.16) ;
%Straight Lines [id:da04381316632330545] 
\draw    (68.03,107.24) -- (68.03,60.4) ;
%Shape: Free Drawing [id:dp8001606872179554] 
\draw  [line width=3] [line join = round][line cap = round] (67.87,107.02) .. controls (68.29,107.02) and (68.29,107.02) .. (67.87,107.02) ;
%Straight Lines [id:da8612850296026056] 
\draw    (101.9,89.89) -- (86.26,45.66) ;
%Shape: Free Drawing [id:dp8357052417274689] 
\draw  [line width=3] [line join = round][line cap = round] (101.74,89.67) .. controls (102.16,89.67) and (102.16,89.67) .. (101.74,89.67) ;
%Straight Lines [id:da9458949653329614] 
\draw    (45.44,52.16) -- (68.03,60.4) ;
%Straight Lines [id:da6500857739397696] 
\draw    (68.03,60.4) -- (86.26,45.66) ;
%Straight Lines [id:da9867695350796011] 
\draw    (45.44,52.16) -- (86.26,45.66) ;
%Straight Lines [id:da3520481628683607] 
\draw    (34.15,99.43) -- (68.03,107.24) ;
%Straight Lines [id:da2551365736268243] 
\draw    (68.03,107.24) -- (101.9,89.89) ;
%Shape: Ellipse [id:dp7960577613453648] 
\draw   (125.33,74.71) .. controls (125.33,49.33) and (145.39,28.75) .. (170.13,28.75) .. controls (194.88,28.75) and (214.94,49.33) .. (214.94,74.71) .. controls (214.94,100.1) and (194.88,120.68) .. (170.13,120.68) .. controls (145.39,120.68) and (125.33,100.1) .. (125.33,74.71) -- cycle ;
%Shape: Arc [id:dp05353246207541251] 
\draw  [draw opacity=0] (214.86,78.19) .. controls (208.32,82.03) and (191.4,84.79) .. (171.57,84.83) .. controls (146.04,84.88) and (125.34,80.42) .. (125.32,74.86) .. controls (125.32,74.81) and (125.33,74.76) .. (125.33,74.71) -- (171.54,74.76) -- cycle ; \draw   (214.86,78.19) .. controls (208.32,82.03) and (191.4,84.79) .. (171.57,84.83) .. controls (146.04,84.88) and (125.34,80.42) .. (125.32,74.86) .. controls (125.32,74.81) and (125.33,74.76) .. (125.33,74.71) ;  
%Shape: Free Drawing [id:dp3796810690775547] 
\draw  [line width=3] [line join = round][line cap = round] (141.69,54.11) .. controls (142.11,54.11) and (142.11,54.11) .. (141.69,54.11) ;
%Shape: Free Drawing [id:dp0765400524807962] 
\draw  [line width=3] [line join = round][line cap = round] (170.35,28.96) .. controls (170.77,28.96) and (170.77,28.96) .. (170.35,28.96) ;
%Shape: Free Drawing [id:dp39693297036814656] 
\draw  [line width=3] [line join = round][line cap = round] (140.82,100.95) .. controls (141.24,100.95) and (141.24,100.95) .. (140.82,100.95) ;
%Shape: Free Drawing [id:dp7688450453979006] 
\draw  [line width=3] [line join = round][line cap = round] (199.01,53.25) .. controls (199.44,53.25) and (199.44,53.25) .. (199.01,53.25) ;
%Shape: Free Drawing [id:dp7056410414713095] 
\draw  [line width=3] [line join = round][line cap = round] (171.22,60.18) .. controls (171.64,60.18) and (171.64,60.18) .. (171.22,60.18) ;
%Shape: Free Drawing [id:dp9754420481116246] 
\draw  [line width=3] [line join = round][line cap = round] (171.22,107.02) .. controls (171.64,107.02) and (171.64,107.02) .. (171.22,107.02) ;
%Shape: Free Drawing [id:dp4767313358392782] 
\draw  [line width=3] [line join = round][line cap = round] (199.01,94.01) .. controls (199.44,94.01) and (199.44,94.01) .. (199.01,94.01) ;
%Straight Lines [id:da3407310346458128] 
\draw    (141.85,54.33) -- (140.98,101.16) ;
%Straight Lines [id:da1353879330597041] 
\draw    (171.38,59.53) -- (171.38,107.24) ;
%Straight Lines [id:da4361923884138299] 
\draw    (199.35,53.03) -- (199.35,94.49) ;
%Straight Lines [id:da7387735133011963] 
\draw    (141.85,54.33) -- (171.38,59.53) ;
%Straight Lines [id:da712581961111932] 
\draw    (171.38,59.53) -- (199.35,53.03) ;
%Straight Lines [id:da001953868293829375] 
\draw    (140.98,101.16) -- (170.51,106.37) ;
%Straight Lines [id:da11163645950731604] 
\draw    (171.38,107.24) -- (199.35,94.49) ;
%Straight Lines [id:da6555402454528994] 
\draw    (141.85,54.33) -- (199.35,53.03) ;
%Shape: Ellipse [id:dp3770087997263356] 
\draw   (320.26,74.73) .. controls (320.49,100.11) and (300.62,120.87) .. (275.87,121.1) .. controls (251.13,121.32) and (230.88,100.93) .. (230.65,75.54) .. controls (230.42,50.15) and (250.29,29.39) .. (275.04,29.17) .. controls (299.78,28.94) and (320.03,49.34) .. (320.26,74.73) -- cycle ;
%Shape: Arc [id:dp3442562727423457] 
\draw  [draw opacity=0] (320.18,79.27) .. controls (313.62,83.09) and (296.7,85.8) .. (276.86,85.78) .. controls (251.33,85.77) and (230.64,81.24) .. (230.65,75.68) .. controls (230.65,75.63) and (230.65,75.59) .. (230.65,75.54) -- (276.87,75.72) -- cycle ; \draw   (320.18,79.27) .. controls (313.62,83.09) and (296.7,85.8) .. (276.86,85.78) .. controls (251.33,85.77) and (230.64,81.24) .. (230.65,75.68) .. controls (230.65,75.63) and (230.65,75.59) .. (230.65,75.54) ;  
%Shape: Free Drawing [id:dp1008973277475873] 
\draw  [line width=3] [line join = round][line cap = round] (297.15,97.27) .. controls (296.73,97.27) and (296.73,97.27) .. (297.15,97.27) ;
%Shape: Free Drawing [id:dp04842821395266106] 
\draw  [line width=3] [line join = round][line cap = round] (275.27,29.17) .. controls (275.41,29.17) and (275.41,29.17) .. (275.27,29.17) ;
%Shape: Free Drawing [id:dp8983727620405508] 
\draw  [line width=3] [line join = round][line cap = round] (308.01,49.47) .. controls (307.59,49.47) and (307.59,49.47) .. (308.01,49.47) ;
%Shape: Free Drawing [id:dp5013829962206084] 
\draw  [line width=3] [line join = round][line cap = round] (256.39,103.71) .. controls (255.97,103.71) and (255.97,103.71) .. (256.39,103.71) ;
%Shape: Free Drawing [id:dp513065445689904] 
\draw  [line width=3] [line join = round][line cap = round] (274.5,89.67) .. controls (274.08,89.67) and (274.08,89.67) .. (274.5,89.67) ;
%Straight Lines [id:da07713366334413219] 
\draw    (307.85,50.12) -- (296.99,97.49) ;
%Straight Lines [id:da9205282308407547] 
\draw    (273.91,42.62) -- (274.34,89.45) ;
%Shape: Free Drawing [id:dp5656290542498836] 
\draw  [line width=3] [line join = round][line cap = round] (274.07,42.84) .. controls (273.65,42.84) and (273.65,42.84) .. (274.07,42.84) ;
%Straight Lines [id:da78783770779622] 
\draw    (240.2,60.27) -- (256.23,104.36) ;
%Shape: Free Drawing [id:dp8908548045101858] 
\draw  [line width=3] [line join = round][line cap = round] (240.36,60.49) .. controls (239.94,60.49) and (239.94,60.49) .. (240.36,60.49) ;
%Straight Lines [id:da8100592885952174] 
\draw    (296.99,97.49) -- (274.34,89.45) ;
%Straight Lines [id:da4821167100639143] 
\draw    (274.34,89.45) -- (256.23,104.36) ;
%Straight Lines [id:da4967530719721085] 
\draw  [dash pattern={on 0.84pt off 2.51pt}]  (307.85,50.12) -- (240.2,60.27) ;
%Straight Lines [id:da5762618294735564] 
\draw    (307.85,50.12) -- (273.91,42.62) ;
%Straight Lines [id:da7733897163170169] 
\draw    (273.91,42.62) -- (240.2,60.27) ;
%Straight Lines [id:da10257094352409846] 
\draw  [dash pattern={on 0.84pt off 2.51pt}]  (101.81,89.28) -- (34.15,99.43) ;
%Straight Lines [id:da5312996618447918] 
\draw    (296.99,97.49) -- (256.23,104.36) ;
%Straight Lines [id:da588132471545693] 
\draw    (140.98,101.16) -- (199.35,94.49) ;
%Shape: Ellipse [id:dp5782218491515789] 
\draw   (376.93,120.79) .. controls (351.99,121.35) and (331.32,101.78) .. (330.76,77.08) .. controls (330.21,52.38) and (349.97,31.9) .. (374.91,31.34) .. controls (399.85,30.78) and (420.52,50.35) .. (421.08,75.05) .. controls (421.63,99.75) and (401.87,120.23) .. (376.93,120.79) -- cycle ;
%Shape: Arc [id:dp5280828205105939] 
\draw  [draw opacity=0] (420.27,81.2) .. controls (413.69,84.99) and (396.76,87.63) .. (376.92,87.53) .. controls (351.4,87.4) and (330.73,82.79) .. (330.76,77.23) .. controls (330.76,77.18) and (330.76,77.13) .. (330.76,77.08) -- (376.97,77.46) -- cycle ; \draw   (420.27,81.2) .. controls (413.69,84.99) and (396.76,87.63) .. (376.92,87.53) .. controls (351.4,87.4) and (330.73,82.79) .. (330.76,77.23) .. controls (330.76,77.18) and (330.76,77.13) .. (330.76,77.08) ;  
%Shape: Free Drawing [id:dp726539452924525] 
\draw  [line width=3] [line join = round][line cap = round] (352.75,98.05) .. controls (352.74,97.63) and (352.74,97.63) .. (352.75,98.05) ;
%Shape: Free Drawing [id:dp51235080581624] 
\draw  [line width=3] [line join = round][line cap = round] (400.76,108.25) .. controls (400.75,107.83) and (400.75,107.83) .. (400.76,108.25) ;
%Shape: Free Drawing [id:dp2404510481137352] 
\draw  [line width=3] [line join = round][line cap = round] (345.75,57.44) .. controls (345.74,57.02) and (345.74,57.02) .. (345.75,57.44) ;
%Shape: Free Drawing [id:dp7247886691125156] 
\draw  [line width=3] [line join = round][line cap = round] (360.06,75.33) .. controls (360.05,74.91) and (360.05,74.91) .. (360.06,75.33) ;
%Straight Lines [id:da5749925483087788] 
\draw    (400.1,108.11) -- (352.53,97.9) ;
%Straight Lines [id:da32701712523345106] 
\draw    (407.16,74.12) -- (360.27,75.17) ;
%Shape: Free Drawing [id:dp7687833329169855] 
\draw  [line width=3] [line join = round][line cap = round] (406.94,74.28) .. controls (406.93,73.86) and (406.93,73.86) .. (406.94,74.28) ;
%Straight Lines [id:da570455076843327] 
\draw    (389.03,40.69) -- (345.1,57.3) ;
%Shape: Free Drawing [id:dp17602067989400672] 
\draw  [line width=3] [line join = round][line cap = round] (388.81,40.85) .. controls (388.81,40.43) and (388.81,40.43) .. (388.81,40.85) ;
%Straight Lines [id:da7515435812013005] 
\draw    (352.53,97.9) -- (360.27,75.17) ;
%Straight Lines [id:da7469106104097574] 
\draw    (360.27,75.17) -- (345.1,57.3) ;
%Straight Lines [id:da4574871111208252] 
\draw  [dash pattern={on 0.84pt off 2.51pt}]  (400.1,108.11) -- (389.03,40.69) ;
%Straight Lines [id:da795204521030511] 
\draw    (400.1,108.11) -- (407.16,74.12) ;
%Straight Lines [id:da018371907170575685] 
\draw    (407.16,74.12) -- (389.03,40.69) ;
%Straight Lines [id:da5906473808323534] 
\draw    (352.53,97.9) -- (345.1,57.3) ;
%Shape: Free Drawing [id:dp6229422940666394] 
\draw  [line width=3] [line join = round][line cap = round] (376.24,31.3) .. controls (376.1,31.31) and (376.1,31.31) .. (376.24,31.3) ;
%Shape: Ellipse [id:dp6721628316366814] 
\draw   (439.01,77.33) .. controls (438.2,52.43) and (457.6,31.6) .. (482.33,30.79) .. controls (507.06,29.99) and (527.77,49.52) .. (528.57,74.42) .. controls (529.38,99.32) and (509.99,120.16) .. (485.26,120.96) .. controls (460.52,121.77) and (439.82,102.23) .. (439.01,77.33) -- cycle ;
%Shape: Arc [id:dp1106935787565253] 
\draw  [draw opacity=0] (528.53,81.19) .. controls (521.97,85) and (505.04,87.68) .. (485.21,87.64) .. controls (459.68,87.59) and (439,83.04) .. (439.01,77.48) .. controls (439.01,77.43) and (439.01,77.38) .. (439.01,77.33) -- (485.23,77.57) -- cycle ; \draw   (528.53,81.19) .. controls (521.97,85) and (505.04,87.68) .. (485.21,87.64) .. controls (459.68,87.59) and (439,83.04) .. (439.01,77.48) .. controls (439.01,77.43) and (439.01,77.38) .. (439.01,77.33) ;  
%Shape: Free Drawing [id:dp2032848020582365] 
\draw  [line width=3] [line join = round][line cap = round] (461.54,52.96) .. controls (461.96,52.95) and (461.96,52.95) .. (461.54,52.96) ;
%Shape: Free Drawing [id:dp8356081236678303] 
\draw  [line width=3] [line join = round][line cap = round] (451.81,101) .. controls (452.23,100.99) and (452.23,100.99) .. (451.81,101) ;
%Shape: Free Drawing [id:dp8021445482613132] 
\draw  [line width=3] [line join = round][line cap = round] (502.15,45.57) .. controls (502.57,45.56) and (502.57,45.56) .. (502.15,45.57) ;
%Shape: Free Drawing [id:dp21193613807974687] 
\draw  [line width=3] [line join = round][line cap = round] (484.37,60.03) .. controls (484.79,60.02) and (484.79,60.02) .. (484.37,60.03) ;
%Straight Lines [id:da2676305471017668] 
\draw    (451.95,100.35) -- (461.69,52.74) ;
%Straight Lines [id:da2806116589709846] 
\draw    (486.05,107.05) -- (484.53,60.24) ;
%Shape: Free Drawing [id:dp35370807063344567] 
\draw  [line width=3] [line join = round][line cap = round] (485.89,106.84) .. controls (486.31,106.83) and (486.31,106.83) .. (485.89,106.84) ;
%Straight Lines [id:da7240851689259469] 
\draw    (519.34,88.62) -- (502.28,44.91) ;
%Shape: Free Drawing [id:dp6089052521574472] 
\draw  [line width=3] [line join = round][line cap = round] (519.18,88.41) .. controls (519.6,88.39) and (519.6,88.39) .. (519.18,88.41) ;
%Straight Lines [id:da7677988885897451] 
\draw    (461.69,52.74) -- (484.53,60.24) ;
%Straight Lines [id:da797195741663782] 
\draw    (484.53,60.24) -- (502.28,44.91) ;
%Straight Lines [id:da8649738831256515] 
\draw  [dash pattern={on 0.84pt off 2.51pt}]  (451.95,100.35) -- (519.34,88.62) ;
%Straight Lines [id:da43464356579202224] 
\draw    (451.95,100.35) -- (486.05,107.05) ;
%Straight Lines [id:da6107254244615774] 
\draw    (486.05,107.05) -- (519.34,88.62) ;
%Straight Lines [id:da24932405949134706] 
\draw    (461.69,52.74) -- (502.28,44.91) ;
%Shape: Free Drawing [id:dp5621798256881139] 
\draw  [line width=3] [line join = round][line cap = round] (483.94,30.44) .. controls (483.8,30.44) and (483.8,30.44) .. (483.94,30.44) ;
%Shape: Ellipse [id:dp7765092851524453] 
\draw   (544.97,78.2) .. controls (544.16,53.3) and (563.56,32.46) .. (588.29,31.66) .. controls (613.02,30.86) and (633.73,50.39) .. (634.53,75.29) .. controls (635.34,100.19) and (615.95,121.03) .. (591.22,121.83) .. controls (566.48,122.63) and (545.78,103.1) .. (544.97,78.2) -- cycle ;
%Shape: Arc [id:dp9880586928501285] 
\draw  [draw opacity=0] (634.49,82.05) .. controls (627.93,85.86) and (611,88.55) .. (591.17,88.51) .. controls (565.64,88.46) and (544.96,83.9) .. (544.97,78.34) .. controls (544.97,78.29) and (544.97,78.25) .. (544.97,78.2) -- (591.19,78.44) -- cycle ; \draw   (634.49,82.05) .. controls (627.93,85.86) and (611,88.55) .. (591.17,88.51) .. controls (565.64,88.46) and (544.96,83.9) .. (544.97,78.34) .. controls (544.97,78.29) and (544.97,78.25) .. (544.97,78.2) ;  
%Shape: Free Drawing [id:dp08992721756907951] 
\draw  [line width=3] [line join = round][line cap = round] (567.5,53.83) .. controls (567.92,53.81) and (567.92,53.81) .. (567.5,53.83) ;
%Shape: Free Drawing [id:dp948098552069393] 
\draw  [line width=3] [line join = round][line cap = round] (557.77,101.87) .. controls (558.19,101.86) and (558.19,101.86) .. (557.77,101.87) ;
%Shape: Free Drawing [id:dp017769083961179044] 
\draw  [line width=3] [line join = round][line cap = round] (608.11,46.44) .. controls (608.53,46.42) and (608.53,46.42) .. (608.11,46.44) ;
%Shape: Free Drawing [id:dp2647158754953516] 
\draw  [line width=3] [line join = round][line cap = round] (590.33,60.9) .. controls (590.75,60.88) and (590.75,60.88) .. (590.33,60.9) ;
%Straight Lines [id:da7258431904921343] 
\draw    (557.91,101.22) -- (567.65,53.61) ;
%Straight Lines [id:da456705994435443] 
\draw    (592.01,107.92) -- (590.49,61.11) ;
%Shape: Free Drawing [id:dp20690810391897374] 
\draw  [line width=3] [line join = round][line cap = round] (591.85,107.71) .. controls (592.27,107.69) and (592.27,107.69) .. (591.85,107.71) ;
%Straight Lines [id:da7910876919520147] 
\draw    (625.3,89.48) -- (608.24,45.78) ;
%Shape: Free Drawing [id:dp12356780554006375] 
\draw  [line width=3] [line join = round][line cap = round] (625.14,89.27) .. controls (625.56,89.26) and (625.56,89.26) .. (625.14,89.27) ;
%Straight Lines [id:da3879236832290698] 
\draw    (567.65,53.61) -- (590.49,61.11) ;
%Straight Lines [id:da81304999291184] 
\draw    (590.49,61.11) -- (608.24,45.78) ;
%Straight Lines [id:da019381489423799136] 
\draw  [dash pattern={on 0.84pt off 2.51pt}]  (557.91,101.22) -- (625.3,89.48) ;
%Straight Lines [id:da47860284817574694] 
\draw    (557.91,101.22) -- (592.01,107.92) ;
%Straight Lines [id:da9446434451519965] 
\draw    (592.01,107.92) -- (625.3,89.48) ;
%Straight Lines [id:da4771097830708313] 
\draw    (567.65,53.61) -- (608.24,45.78) ;
%Shape: Free Drawing [id:dp8484274534370522] 
\draw  [line width=3] [line join = round][line cap = round] (589.9,31.3) .. controls (589.76,31.31) and (589.76,31.31) .. (589.9,31.3) ;
%Shape: Parallelogram [id:dp6490115575163433] 
\draw   (198.3,210.32) -- (114.4,210.76) -- (132.72,149.57) -- (216.62,149.13) -- cycle ;
%Shape: Ellipse [id:dp924376283845127] 
\draw   (130.51,179.95) .. controls (130.51,168.9) and (146.18,159.95) .. (165.51,159.95) .. controls (184.84,159.95) and (200.51,168.9) .. (200.51,179.95) .. controls (200.51,190.99) and (184.84,199.95) .. (165.51,199.95) .. controls (146.18,199.95) and (130.51,190.99) .. (130.51,179.95) -- cycle ;
%Shape: Ellipse [id:dp42251202030348467] 
\draw   (147.27,179.95) .. controls (147.27,173.75) and (155.43,168.72) .. (165.51,168.72) .. controls (175.59,168.72) and (183.76,173.75) .. (183.76,179.95) .. controls (183.76,186.15) and (175.59,191.17) .. (165.51,191.17) .. controls (155.43,191.17) and (147.27,186.15) .. (147.27,179.95) -- cycle ;
%Shape: Free Drawing [id:dp6064956955977332] 
\draw  [line width=3] [line join = round][line cap = round] (141.29,165.38) .. controls (141.71,165.38) and (141.71,165.38) .. (141.29,165.38) ;
%Shape: Free Drawing [id:dp5669630947286617] 
\draw  [line width=3] [line join = round][line cap = round] (185.29,163.38) .. controls (185.71,163.38) and (185.71,163.38) .. (185.29,163.38) ;
%Shape: Free Drawing [id:dp12686517566955247] 
\draw  [line width=3] [line join = round][line cap = round] (152.29,172.38) .. controls (152.71,172.38) and (152.71,172.38) .. (152.29,172.38) ;
%Shape: Free Drawing [id:dp9497340969026213] 
\draw  [line width=3] [line join = round][line cap = round] (176.29,170.38) .. controls (176.71,170.38) and (176.71,170.38) .. (176.29,170.38) ;
%Shape: Free Drawing [id:dp8152821384396962] 
\draw  [line width=3] [line join = round][line cap = round] (165.29,191.38) .. controls (165.71,191.38) and (165.71,191.38) .. (165.29,191.38) ;
%Shape: Free Drawing [id:dp2669308083014139] 
\draw  [line width=3] [line join = round][line cap = round] (165.29,199.38) .. controls (165.71,199.38) and (165.71,199.38) .. (165.29,199.38) ;
%Straight Lines [id:da12488032720240594] 
\draw    (40,136) -- (301.5,135) ;
%Straight Lines [id:da5048261238655647] 
\draw    (40,136) -- (40.5,126) ;
%Straight Lines [id:da9922630759386604] 
\draw    (301.5,135) -- (301.5,124) ;
%Straight Lines [id:da521103124844716] 
\draw    (170.75,146.5) -- (170.75,135.5) ;
%Shape: Parallelogram [id:dp6933532192044773] 
\draw   (303.67,210.6) -- (219.77,211.04) -- (238.09,149.85) -- (321.99,149.41) -- cycle ;
%Shape: Ellipse [id:dp42059543144699707] 
\draw   (236.27,183.95) .. controls (236.27,177.75) and (244.43,172.72) .. (254.51,172.72) .. controls (264.59,172.72) and (272.76,177.75) .. (272.76,183.95) .. controls (272.76,190.15) and (264.59,195.17) .. (254.51,195.17) .. controls (244.43,195.17) and (236.27,190.15) .. (236.27,183.95) -- cycle ;
%Shape: Free Drawing [id:dp8612572871382522] 
\draw  [line width=3] [line join = round][line cap = round] (289.29,166.38) .. controls (289.71,166.38) and (289.71,166.38) .. (289.29,166.38) ;
%Shape: Free Drawing [id:dp05153450132361037] 
\draw  [line width=3] [line join = round][line cap = round] (289.29,180.38) .. controls (289.71,180.38) and (289.71,180.38) .. (289.29,180.38) ;
%Shape: Free Drawing [id:dp3641969619610368] 
\draw  [line width=3] [line join = round][line cap = round] (241.29,176.38) .. controls (241.71,176.38) and (241.71,176.38) .. (241.29,176.38) ;
%Shape: Free Drawing [id:dp7776258695504006] 
\draw  [line width=3] [line join = round][line cap = round] (272.29,182.38) .. controls (272.71,182.38) and (272.71,182.38) .. (272.29,182.38) ;
%Shape: Free Drawing [id:dp7858943134168448] 
\draw  [line width=3] [line join = round][line cap = round] (252.29,195.38) .. controls (252.71,195.38) and (252.71,195.38) .. (252.29,195.38) ;
%Shape: Free Drawing [id:dp6565382062413496] 
\draw  [line width=3] [line join = round][line cap = round] (289.29,196.38) .. controls (289.71,196.38) and (289.71,196.38) .. (289.29,196.38) ;
%Straight Lines [id:da7314297328357724] 
\draw    (289.5,157) -- (289.5,206) ;
%Shape: Parallelogram [id:dp6588634338988113] 
\draw   (405.54,209.32) -- (321.65,209.76) -- (339.97,148.57) -- (423.86,148.13) -- cycle ;
%Shape: Ellipse [id:dp9615061740847799] 
\draw   (336.27,178.95) .. controls (336.27,172.75) and (344.43,167.72) .. (354.51,167.72) .. controls (364.59,167.72) and (372.76,172.75) .. (372.76,178.95) .. controls (372.76,185.15) and (364.59,190.17) .. (354.51,190.17) .. controls (344.43,190.17) and (336.27,185.15) .. (336.27,178.95) -- cycle ;
%Shape: Free Drawing [id:dp9182246144098212] 
\draw  [line width=3] [line join = round][line cap = round] (341.29,171.38) .. controls (341.71,171.38) and (341.71,171.38) .. (341.29,171.38) ;
%Shape: Free Drawing [id:dp7147607491982098] 
\draw  [line width=3] [line join = round][line cap = round] (365.29,169.38) .. controls (365.71,169.38) and (365.71,169.38) .. (365.29,169.38) ;
%Shape: Free Drawing [id:dp9734309487259931] 
\draw  [line width=3] [line join = round][line cap = round] (352.29,190.38) .. controls (352.71,190.38) and (352.71,190.38) .. (352.29,190.38) ;
%Shape: Ellipse [id:dp9568409944415113] 
\draw   (375.27,178.95) .. controls (375.27,172.75) and (383.43,167.72) .. (393.51,167.72) .. controls (403.59,167.72) and (411.76,172.75) .. (411.76,178.95) .. controls (411.76,185.15) and (403.59,190.17) .. (393.51,190.17) .. controls (383.43,190.17) and (375.27,185.15) .. (375.27,178.95) -- cycle ;
%Shape: Free Drawing [id:dp447987376217339] 
\draw  [line width=3] [line join = round][line cap = round] (384.29,169.38) .. controls (384.71,169.38) and (384.71,169.38) .. (384.29,169.38) ;
%Shape: Free Drawing [id:dp07495779796756286] 
\draw  [line width=3] [line join = round][line cap = round] (407.29,171.38) .. controls (407.71,171.38) and (407.71,171.38) .. (407.29,171.38) ;
%Shape: Free Drawing [id:dp7297149725810261] 
\draw  [line width=3] [line join = round][line cap = round] (391.29,190.38) .. controls (391.71,190.38) and (391.71,190.38) .. (391.29,190.38) ;
%Straight Lines [id:da7565673447273975] 
\draw    (20,108) -- (20.47,139) ;
\draw [shift={(20.5,141)}, rotate = 269.13] [color={rgb, 255:red, 0; green, 0; blue, 0 }  ][line width=0.75]    (10.93,-3.29) .. controls (6.95,-1.4) and (3.31,-0.3) .. (0,0) .. controls (3.31,0.3) and (6.95,1.4) .. (10.93,3.29)   ;
%Straight Lines [id:da9826147774824956] 
\draw    (141,165) -- (152,172) ;
%Straight Lines [id:da817793797685763] 
\draw    (177,170) -- (186,162) ;
%Straight Lines [id:da4971984995037364] 
\draw    (165.51,191.17) -- (165.51,199.95) ;
%Straight Lines [id:da25068951002357887] 
\draw    (240,175.5) -- (290,166.5) ;
%Straight Lines [id:da8924231057172667] 
\draw    (272.76,183.95) -- (289.5,181.5) ;
%Straight Lines [id:da8334171735249847] 
\draw    (254.51,195.17) -- (290,196.33) ;
%Straight Lines [id:da18883308825983125] 
\draw    (365,168.5) -- (384,169.5) ;
%Straight Lines [id:da47729303223146635] 
\draw    (340,170.5) -- (408,171.5) ;
%Straight Lines [id:da9487008464381899] 
\draw    (354.51,190.17) -- (393.51,190.17) ;
%Shape: Parallelogram [id:dp06603241590942965] 
\draw   (516.3,208.32) -- (432.4,208.76) -- (450.72,147.57) -- (534.62,147.13) -- cycle ;
%Shape: Ellipse [id:dp40276996826997413] 
\draw   (448.51,177.95) .. controls (448.51,166.9) and (464.18,157.95) .. (483.51,157.95) .. controls (502.84,157.95) and (518.51,166.9) .. (518.51,177.95) .. controls (518.51,188.99) and (502.84,197.95) .. (483.51,197.95) .. controls (464.18,197.95) and (448.51,188.99) .. (448.51,177.95) -- cycle ;
%Shape: Ellipse [id:dp9950985355582287] 
\draw   (465.27,177.95) .. controls (465.27,171.75) and (473.43,166.72) .. (483.51,166.72) .. controls (493.59,166.72) and (501.76,171.75) .. (501.76,177.95) .. controls (501.76,184.15) and (493.59,189.17) .. (483.51,189.17) .. controls (473.43,189.17) and (465.27,184.15) .. (465.27,177.95) -- cycle ;
%Shape: Free Drawing [id:dp04182095223519933] 
\draw  [line width=3] [line join = round][line cap = round] (459.29,163.38) .. controls (459.71,163.38) and (459.71,163.38) .. (459.29,163.38) ;
%Shape: Free Drawing [id:dp8836077555286301] 
\draw  [line width=3] [line join = round][line cap = round] (503.29,161.38) .. controls (503.71,161.38) and (503.71,161.38) .. (503.29,161.38) ;
%Shape: Free Drawing [id:dp8758575956901025] 
\draw  [line width=3] [line join = round][line cap = round] (470.29,170.38) .. controls (470.71,170.38) and (470.71,170.38) .. (470.29,170.38) ;
%Shape: Free Drawing [id:dp11720372552596259] 
\draw  [line width=3] [line join = round][line cap = round] (494.29,168.38) .. controls (494.71,168.38) and (494.71,168.38) .. (494.29,168.38) ;
%Shape: Free Drawing [id:dp0007948866355944606] 
\draw  [line width=3] [line join = round][line cap = round] (483.29,189.38) .. controls (483.71,189.38) and (483.71,189.38) .. (483.29,189.38) ;
%Shape: Free Drawing [id:dp1451536692288944] 
\draw  [line width=3] [line join = round][line cap = round] (483.29,197.38) .. controls (483.71,197.38) and (483.71,197.38) .. (483.29,197.38) ;
%Straight Lines [id:da5863511203930643] 
\draw    (459,163) -- (470,170) ;
%Straight Lines [id:da7359479737144645] 
\draw    (495,168) -- (504,160) ;
%Straight Lines [id:da147673256553293] 
\draw    (483.51,189.17) -- (483.51,197.95) ;
%Shape: Parallelogram [id:dp5702987316302544] 
\draw   (624.3,208.32) -- (540.4,208.76) -- (558.72,147.57) -- (642.62,147.13) -- cycle ;
%Shape: Ellipse [id:dp6570112371935106] 
\draw   (556.51,177.95) .. controls (556.51,166.9) and (572.18,157.95) .. (591.51,157.95) .. controls (610.84,157.95) and (626.51,166.9) .. (626.51,177.95) .. controls (626.51,188.99) and (610.84,197.95) .. (591.51,197.95) .. controls (572.18,197.95) and (556.51,188.99) .. (556.51,177.95) -- cycle ;
%Shape: Ellipse [id:dp7586400490124368] 
\draw   (573.27,177.95) .. controls (573.27,171.75) and (581.43,166.72) .. (591.51,166.72) .. controls (601.59,166.72) and (609.76,171.75) .. (609.76,177.95) .. controls (609.76,184.15) and (601.59,189.17) .. (591.51,189.17) .. controls (581.43,189.17) and (573.27,184.15) .. (573.27,177.95) -- cycle ;
%Shape: Free Drawing [id:dp5610234614447795] 
\draw  [line width=3] [line join = round][line cap = round] (567.29,163.38) .. controls (567.71,163.38) and (567.71,163.38) .. (567.29,163.38) ;
%Shape: Free Drawing [id:dp47039368271838544] 
\draw  [line width=3] [line join = round][line cap = round] (611.29,161.38) .. controls (611.71,161.38) and (611.71,161.38) .. (611.29,161.38) ;
%Shape: Free Drawing [id:dp6514829775322744] 
\draw  [line width=3] [line join = round][line cap = round] (578.29,170.38) .. controls (578.71,170.38) and (578.71,170.38) .. (578.29,170.38) ;
%Shape: Free Drawing [id:dp510409278806826] 
\draw  [line width=3] [line join = round][line cap = round] (602.29,168.38) .. controls (602.71,168.38) and (602.71,168.38) .. (602.29,168.38) ;
%Shape: Free Drawing [id:dp7468746431960029] 
\draw  [line width=3] [line join = round][line cap = round] (591.29,189.38) .. controls (591.71,189.38) and (591.71,189.38) .. (591.29,189.38) ;
%Shape: Free Drawing [id:dp9847075648373711] 
\draw  [line width=3] [line join = round][line cap = round] (591.29,197.38) .. controls (591.71,197.38) and (591.71,197.38) .. (591.29,197.38) ;
%Straight Lines [id:da6469217421152877] 
\draw    (567,163) -- (578,170) ;
%Straight Lines [id:da9844901140582316] 
\draw    (603,168) -- (612,160) ;
%Straight Lines [id:da8308072747691384] 
\draw    (591.51,189.17) -- (591.51,197.95) ;
%Curve Lines [id:da6253199262169207] 
\draw    (301,150) .. controls (312.9,136.98) and (327.54,132.88) .. (326.24,101.66) ;
\draw [shift={(326.14,99.71)}, rotate = 86.53] [color={rgb, 255:red, 0; green, 0; blue, 0 }  ][line width=0.75]    (10.93,-3.29) .. controls (6.95,-1.4) and (3.31,-0.3) .. (0,0) .. controls (3.31,0.3) and (6.95,1.4) .. (10.93,3.29)   ;
%Curve Lines [id:da7960309771506868] 
\draw    (301,150) .. controls (447.39,122.42) and (432.99,121.73) .. (430.32,97.62) ;
\draw [shift={(430.14,95.71)}, rotate = 85.6] [color={rgb, 255:red, 0; green, 0; blue, 0 }  ][line width=0.75]    (10.93,-3.29) .. controls (6.95,-1.4) and (3.31,-0.3) .. (0,0) .. controls (3.31,0.3) and (6.95,1.4) .. (10.93,3.29)   ;

% Text Node
\draw (67.24,16.4) node [anchor=north west][inner sep=0.75pt]  [font=\scriptsize]  {$I$};
% Text Node
\draw (35.1,102.1) node [anchor=north west][inner sep=0.75pt]  [font=\tiny]  {${\textstyle 1}$};
% Text Node
\draw (68.97,63.07) node [anchor=north west][inner sep=0.75pt]  [font=\tiny]  {$2$};
% Text Node
\draw (96.77,93.43) node [anchor=north west][inner sep=0.75pt]  [font=\tiny]  {$3$};
% Text Node
\draw (46.39,54.83) node [anchor=north west][inner sep=0.75pt]  [font=\tiny]  {$4$};
% Text Node
\draw (59.42,109.44) node [anchor=north west][inner sep=0.75pt]  [font=\tiny]  {$ \begin{array}{l}
5\\
\end{array}$};
% Text Node
\draw (81.13,51.8) node [anchor=north west][inner sep=0.75pt]  [font=\tiny]  {$6$};
% Text Node
\draw (170.59,15.54) node [anchor=north west][inner sep=0.75pt]  [font=\scriptsize]  {$I$};
% Text Node
\draw (141.06,92.56) node [anchor=north west][inner sep=0.75pt]  [font=\tiny]  {${\textstyle 1}$};
% Text Node
\draw (172.33,63.07) node [anchor=north west][inner sep=0.75pt]  [font=\tiny]  {$2$};
% Text Node
\draw (190.57,87.36) node [anchor=north west][inner sep=0.75pt]  [font=\tiny]  {$3$};
% Text Node
\draw (173.38,109.64) node [anchor=north west][inner sep=0.75pt]  [font=\tiny]  {$ \begin{array}{l}
5\\
\end{array}$};
% Text Node
\draw (190.57,58.73) node [anchor=north west][inner sep=0.75pt]  [font=\tiny]  {$6$};
% Text Node
\draw (142.8,57) node [anchor=north west][inner sep=0.75pt]  [font=\tiny]  {$4$};
% Text Node
\draw (275.68,16.4) node [anchor=north west][inner sep=0.75pt]  [font=\scriptsize]  {$I$};
% Text Node
\draw (245.28,60.9) node [anchor=north west][inner sep=0.75pt]  [font=\tiny]  {$4$};
% Text Node
\draw (257.18,107.03) node [anchor=north west][inner sep=0.75pt]  [font=\tiny]  {${\textstyle 1}$};
% Text Node
\draw (274.86,45.29) node [anchor=north west][inner sep=0.75pt]  [font=\tiny]  {$2$};
% Text Node
\draw (276.34,92.85) node [anchor=north west][inner sep=0.75pt]  [font=\tiny]  {$ \begin{array}{l}
5\\
\end{array}$};
% Text Node
\draw (298.26,53.53) node [anchor=north west][inner sep=0.75pt]  [font=\tiny]  {$6$};
% Text Node
\draw (292.18,101.23) node [anchor=north west][inner sep=0.75pt]  [font=\tiny]  {$3$};
% Text Node
\draw (382.11,27.61) node [anchor=north west][inner sep=0.75pt]  [font=\scriptsize,rotate=-179.14]  {$I$};
% Text Node
\draw (350.11,56.73) node [anchor=north west][inner sep=0.75pt]  [font=\tiny]  {${\textstyle 1}$};
% Text Node
\draw (361.22,77.11) node [anchor=north west][inner sep=0.75pt]  [font=\tiny]  {$ \begin{array}{l}
5\\
\end{array}$};
% Text Node
\draw (353.48,100.57) node [anchor=north west][inner sep=0.75pt]  [font=\tiny]  {$3$};
% Text Node
\draw (393.8,40.09) node [anchor=north west][inner sep=0.75pt]  [font=\tiny]  {$4$};
% Text Node
\draw (408.61,66.11) node [anchor=north west][inner sep=0.75pt]  [font=\tiny]  {$2$};
% Text Node
\draw (402.49,99.5) node [anchor=north west][inner sep=0.75pt]  [font=\tiny]  {$6$};
% Text Node
\draw (497.51,51.95) node [anchor=north west][inner sep=0.75pt]  [font=\tiny,rotate=-1.41]  {${\textstyle 1}$};
% Text Node
\draw (476.17,63.01) node [anchor=north west][inner sep=0.75pt]  [font=\tiny,rotate=-1.3]  {$ \begin{array}{l}
5\\
\end{array}$};
% Text Node
\draw (452.76,54.04) node [anchor=north west][inner sep=0.75pt]  [font=\tiny,rotate=-1.09]  {$3$};
% Text Node
\draw (515.45,93.77) node [anchor=north west][inner sep=0.75pt]  [font=\tiny,rotate=-0.52]  {$4$};
% Text Node
\draw (487.02,109.71) node [anchor=north west][inner sep=0.75pt]  [font=\tiny,rotate=-359.53]  {$2$};
% Text Node
\draw (452.83,103.04) node [anchor=north west][inner sep=0.75pt]  [font=\tiny,rotate=-1.34]  {$6$};
% Text Node
\draw (491.54,27.61) node [anchor=north west][inner sep=0.75pt]  [font=\scriptsize,rotate=-179.14]  {$I$};
% Text Node
\draw (582.63,64.09) node [anchor=north west][inner sep=0.75pt]  [font=\tiny,rotate=-1.41]  {${\textstyle 1}$};
% Text Node
\draw (603.84,51.73) node [anchor=north west][inner sep=0.75pt]  [font=\tiny,rotate=-1.3]  {$ \begin{array}{l}
5\\
\end{array}$};
% Text Node
\draw (558.72,54.91) node [anchor=north west][inner sep=0.75pt]  [font=\tiny,rotate=-1.09]  {$3$};
% Text Node
\draw (586.67,111.98) node [anchor=north west][inner sep=0.75pt]  [font=\tiny,rotate=-0.52]  {$4$};
% Text Node
\draw (619.91,94.1) node [anchor=north west][inner sep=0.75pt]  [font=\tiny,rotate=-359.53]  {$2$};
% Text Node
\draw (558.79,103.91) node [anchor=north west][inner sep=0.75pt]  [font=\tiny,rotate=-1.34]  {$6$};
% Text Node
\draw (597.5,28.48) node [anchor=north west][inner sep=0.75pt]  [font=\scriptsize,rotate=-179.14]  {$I$};
% Text Node
\draw (152.1,175.1) node [anchor=north west][inner sep=0.75pt]  [font=\tiny]  {${\textstyle 1}$};
% Text Node
\draw (167.51,183.35) node [anchor=north west][inner sep=0.75pt]  [font=\tiny]  {$ \begin{array}{l}
5\\
\end{array}$};
% Text Node
\draw (167.51,172.12) node [anchor=north west][inner sep=0.75pt]  [font=\tiny]  {$3$};
% Text Node
\draw (139.39,155.83) node [anchor=north west][inner sep=0.75pt]  [font=\tiny]  {$4$};
% Text Node
\draw (188.13,155.8) node [anchor=north west][inner sep=0.75pt]  [font=\tiny]  {$6$};
% Text Node
\draw (167.51,203.35) node [anchor=north west][inner sep=0.75pt]  [font=\tiny]  {$2$};
% Text Node
\draw (241.1,179.1) node [anchor=north west][inner sep=0.75pt]  [font=\tiny]  {${\textstyle 1}$};
% Text Node
\draw (250.51,185.35) node [anchor=north west][inner sep=0.75pt]  [font=\tiny]  {$ \begin{array}{l}
5\\
\end{array}$};
% Text Node
\draw (263.09,179.25) node [anchor=north west][inner sep=0.75pt]  [font=\tiny]  {$3$};
% Text Node
\draw (291.5,160.4) node [anchor=north west][inner sep=0.75pt]  [font=\tiny]  {$4$};
% Text Node
\draw (291.13,173.8) node [anchor=north west][inner sep=0.75pt]  [font=\tiny]  {$6$};
% Text Node
\draw (291.51,190.35) node [anchor=north west][inner sep=0.75pt]  [font=\tiny]  {$2$};
% Text Node
\draw (341.1,174.1) node [anchor=north west][inner sep=0.75pt]  [font=\tiny]  {${\textstyle 1}$};
% Text Node
\draw (350.51,180.35) node [anchor=north west][inner sep=0.75pt]  [font=\tiny]  {$ \begin{array}{l}
5\\
\end{array}$};
% Text Node
\draw (356.51,171.12) node [anchor=north west][inner sep=0.75pt]  [font=\tiny]  {$3$};
% Text Node
\draw (5,114.4) node [anchor=north west][inner sep=0.75pt]  [font=\scriptsize]  {$S_{I}$};
% Text Node
\draw (402.5,175.4) node [anchor=north west][inner sep=0.75pt]  [font=\tiny]  {$4$};
% Text Node
\draw (385.13,170.8) node [anchor=north west][inner sep=0.75pt]  [font=\tiny]  {$6$};
% Text Node
\draw (388.51,180.35) node [anchor=north west][inner sep=0.75pt]  [font=\tiny]  {$2$};
% Text Node
\draw (508.1,155.1) node [anchor=north west][inner sep=0.75pt]  [font=\tiny]  {${\textstyle 1}$};
% Text Node
\draw (473.51,199.35) node [anchor=north west][inner sep=0.75pt]  [font=\tiny]  {$ \begin{array}{l}
5\\
\end{array}$};
% Text Node
\draw (455.51,155.12) node [anchor=north west][inner sep=0.75pt]  [font=\tiny]  {$3$};
% Text Node
\draw (485.51,170.12) node [anchor=north west][inner sep=0.75pt]  [font=\tiny]  {$4$};
% Text Node
\draw (472,173.4) node [anchor=north west][inner sep=0.75pt]  [font=\tiny]  {$6$};
% Text Node
\draw (485.51,181.35) node [anchor=north west][inner sep=0.75pt]  [font=\tiny]  {$2$};
% Text Node
\draw (599.51,198.35) node [anchor=north west][inner sep=0.75pt]  [font=\tiny]  {${\textstyle 1}$};
% Text Node
\draw (615.51,154.35) node [anchor=north west][inner sep=0.75pt]  [font=\tiny]  {$ \begin{array}{l}
5\\
\end{array}$};
% Text Node
\draw (563.51,155.12) node [anchor=north west][inner sep=0.75pt]  [font=\tiny]  {$3$};
% Text Node
\draw (593.51,181.35) node [anchor=north west][inner sep=0.75pt]  [font=\tiny]  {$4$};
% Text Node
\draw (580,173.4) node [anchor=north west][inner sep=0.75pt]  [font=\tiny]  {$6$};
% Text Node
\draw (598.51,172.35) node [anchor=north west][inner sep=0.75pt]  [font=\tiny]  {$2$};
% Text Node
\draw (150,217.4) node [anchor=north west][inner sep=0.75pt]    {$( 1)$};
% Text Node
\draw (254,217.4) node [anchor=north west][inner sep=0.75pt]    {$( 2)$};
% Text Node
\draw (355,216.4) node [anchor=north west][inner sep=0.75pt]    {$( 3)$};
% Text Node
\draw (469,214.4) node [anchor=north west][inner sep=0.75pt]    {$( 4)$};
% Text Node
\draw (577,212.4) node [anchor=north west][inner sep=0.75pt]    {$( 5)$};

\end{tikzpicture}

\end{center}
\caption{All Planar Configurations Corresponding to Elements of \([D \cap Q_\gamma]\)}
\end{figure}

\begin{lemma}
The stereographic projection \(S_I\) has the property that the intersection of the line passing through \(p\) and \(I\), \(\ell\), and the plane \(S_I(T)\) with \(I \in T\), a copy of \(S^2\), is exactly the point of intersection of \(S_I(\ell_1), S_I(\ell_2), S_I(\ell_3)\).
\end{lemma}
\par \textit{Proof.} Since the points \(p_1,p_2, \ldots, p_6\) intersect at \(p\), the trace of the line passing through the center of \(T\) and a point on any \(l_i\) consists of three circles on \(T\) that intersect at two points. At least one of these points must not be \(I\), call this point \(J\), so taking \(J\) under \(S_I\) is the intersection between \(S_I(\ell_1), S_I(\ell_2), S_I(\ell_3)\). However, if both of these points were not I, then the set of images of lines under \(S_I\) would intersect at two points, which is impossible since the lines are not all identical. Therefore, the second point on the intersection of these three circles must be \(I\), such that \(I\), \(J\), and the intersection between \(S_I(\ell_1), S_I(\ell_2), S_I(\ell_3)\) are colinear. Lastly, drawing a line from \(J\) through points on the circle spanned by \(I\) and pairs \((p_1,p_4)\), \((p_2,p_5)\), \((p_3,p_6)\) shows that \(I,J,\) and \(p\) are colinear. Therefore, \(p, I\), and the intersection between \(S_I(\ell_1), S_I(\ell_2), S_I(\ell_3)\) are colinear. \(\boxed{}\)
\begin{definition}
Let \(r_{\ell,\epsilon}(x)\) denote the \(\epsilon\)-rotation of a configuration \(x (\in S_I([D' \cap Q_\gamma])) \subset P\) about \(\ell \in P\), for some plane \(P\). Then, let  \(C_{\ell, \epsilon, i} = \{\pi_i(r_{\ell,\epsilon}(x)) | |\epsilon'| \leq \epsilon\}\), where \(\pi_i(x)\) is the \(i^{\text{th}}\) coordinate of \(x\) corresponding to that \(x\)'s \(p_i\). If \(\forall \ell\) \(\exists \{\epsilon_j\} \rightarrow 0\) such that \(\forall j\), connecting some \(v_i \in C_{\ell, \epsilon_j, i}\) in cyclic order of \(i \in [1,6]\) yields a trefoil, we call \(x\) a \textit{good} configuration. Otherwise, \(x\) is a \textit{bad} configuration.
\end{definition}
\begin{lemma}
 For any configuration \(x\) of types (1), (4), or (5), if \(\epsilon\) is the amount by which the outer circle is rotated with respect to the corresponding \(x' \in [D \cap Q_\gamma]\), an arbitrary rotation \(\epsilon' << \epsilon\) of a single point on the outer circle along it will take a bad configuration to a good one. Furthermore, if \(x\) is good, then there exists an interval \(I \subset [-\epsilon',\epsilon']\) for which such a rotation by any value in this interval yields a good configuration.
\end{lemma}

\textit{Proof.} Suppose \(x\) is a bad configuration. Since \(x\) lies within two nested circles, Lemma 3.7 guarantees that the intersection point of lines \(\ell_{1,2,3}\) must lie within the inner circle. If there are three colinear points in \(x\), the intersection point of all lines must lie on a colinear line, and there must be four colinear points. However, since all combinations of four points either have three on an inner or outer circle, or four on a non-degenerate circle corresponding to the side of the triangular prism in \(M\), we have a contradiction. Thus, there are at most two points of the six on the same line. Picking \(\ell\) passing through at most one point of \(x\) yields a trefoil for small \(\epsilon'\), so one of \(p_{1,4}\), \(p_{2,5}\), or \(p_{3,6}\) lies on \(\ell\) \(\in P \supset x\). Without loss of generality, let it be \(p_{3,6}\). Then, according to Figure 6, where approximately-orthogonal heights to \(P\) are labelled with \(\epsilon_i\) at \(p_i\) and the length of segments up until the point of intersection of lines are labelled with \(a_i, b_i, l_i\), we need relations for the ordering of segments to form a trefoil: \(\frac{l_1}{l_1+l_2}\epsilon_5 > \frac{l_3}{l_3+l_4} \epsilon_2\), \(\frac{a_1}{a_1+a_2} \epsilon_1 + \frac{a_2}{a_1+a_2} \epsilon_2 > \frac{a_3}{a_3+a_4} \epsilon_4 + \frac{a_4}{a_3+a_4} \epsilon_5\), \(\frac{b_3}{b_3+b_4} \epsilon_4 > \frac{b_2}{b_1+b_2} \epsilon_1\), for \(\epsilon_3, \epsilon_6 = 0\). If all of the inequality signs are reversed, we have a trefoil of the opposite handedness such that the only condition that must be satisfied is \(\frac{a_1}{a_1+a_2} \frac{b_3}{b_2} \frac{b_1+b_2}{b_3+b_4} \epsilon_4 + \frac{a_2}{a_1+a_2} \frac{l_1}{l_3} \frac{l_3+l_4}{l_1+l_2} \epsilon_5 >/< \frac{a_1}{a_1+a_2} \epsilon_1 + \frac{a_2}{a_1+a_2} \epsilon_2 >/< \frac{a_3}{a_3+a_4} \epsilon_4 + \frac{a_4}{a_3+a_4} \epsilon_5\), with \(>/<\) representing a joint sign-change in the inequalities. By choosing \(\epsilon_1, \epsilon_2\) appropriately, this condition reduces to \(\neg (\frac{\tilde{a_4}}{\tilde{a_2}} = \frac{\tilde{l_1}}{\tilde{l_3}} \wedge \frac{\tilde{a_3}}{\tilde{a_1}} = \frac{\tilde{b_3}}{\tilde{b_2}})\), where the tilde represents the fraction of the larger line segment that the segment corresponding to the length represents. For instance, \(\tilde{a_4} = \frac{a_4}{a_3 + a_4}\).

\begin{center}
\begin{figure}
\tikzset{every picture/.style={line width=0.75pt}} %set default line width to 0.75pt     

\begin{tikzpicture}[x=0.75pt,y=0.75pt,yscale=-1,xscale=1]
%Shape: Circle [id:dp3913830377348415] 
\draw   (143,147.71) .. controls (143,79.94) and (197.94,25) .. (265.71,25) .. controls (333.49,25) and (388.43,79.94) .. (388.43,147.71) .. controls (388.43,215.49) and (333.49,270.43) .. (265.71,270.43) .. controls (197.94,270.43) and (143,215.49) .. (143,147.71) -- cycle ;
%Shape: Circle [id:dp5296441336722646] 
\draw   (177.5,147.71) .. controls (177.5,98.99) and (216.99,59.5) .. (265.71,59.5) .. controls (314.43,59.5) and (353.93,98.99) .. (353.93,147.71) .. controls (353.93,196.43) and (314.43,235.93) .. (265.71,235.93) .. controls (216.99,235.93) and (177.5,196.43) .. (177.5,147.71) -- cycle ;
%Straight Lines [id:da6434385048773301] 
\draw  [dash pattern={on 0.84pt off 2.51pt}]  (173,68) -- (198.43,90.71) ;
%Straight Lines [id:da8273122475300958] 
\draw  [dash pattern={on 0.84pt off 2.51pt}]  (265.71,235.93) -- (265.71,270.43) ;
%Straight Lines [id:da9330545258739855] 
\draw  [dash pattern={on 0.84pt off 2.51pt}]  (356.43,65.71) -- (332.43,89.71) ;
%Shape: Free Drawing [id:dp7580359168980584] 
\draw  [line width=3] [line join = round][line cap = round] (172.43,67.71) .. controls (172.43,67.71) and (172.43,67.71) .. (172.43,67.71) ;
%Shape: Free Drawing [id:dp3544803835501009] 
\draw  [line width=3] [line join = round][line cap = round] (198.43,91.71) .. controls (198.43,91.71) and (198.43,91.71) .. (198.43,91.71) ;
%Shape: Free Drawing [id:dp5434901428386563] 
\draw  [line width=3] [line join = round][line cap = round] (332.43,88.71) .. controls (332.43,88.71) and (332.43,88.71) .. (332.43,88.71) ;
%Shape: Free Drawing [id:dp27681594129305864] 
\draw  [line width=3] [line join = round][line cap = round] (356.43,65.71) .. controls (356.43,65.71) and (356.43,65.71) .. (356.43,65.71) ;
%Shape: Free Drawing [id:dp9637660060914786] 
\draw  [line width=3] [line join = round][line cap = round] (265.43,235.71) .. controls (265.43,235.71) and (265.43,235.71) .. (265.43,235.71) ;
%Shape: Free Drawing [id:dp6111714151946457] 
\draw  [line width=3] [line join = round][line cap = round] (265.43,269.71) .. controls (265.43,269.71) and (265.43,269.71) .. (265.43,269.71) ;
%Straight Lines [id:da09644269107752601] 
\draw    (198.43,90.71) -- (265.71,270.43) ;
%Straight Lines [id:da2304131884026377] 
\draw    (299.43,179.29) -- (265.71,270.43) ;
%Straight Lines [id:da8273803562692137] 
\draw    (173,68) -- (332.43,89.71) ;
%Straight Lines [id:da9867410761264255] 
\draw    (173,68) -- (214.43,145.29) ;
%Straight Lines [id:da8291210838397638] 
\draw    (356.43,65.71) -- (265.71,235.93) ;
%Straight Lines [id:da7992450900687469] 
\draw    (275.43,77.29) -- (356.43,65.71) ;
%Straight Lines [id:da7960789215536879] 
\draw    (332.43,89.71) -- (316.43,132.29) ;
%Straight Lines [id:da3415932219455009] 
\draw    (237.43,184.29) -- (265.71,235.93) ;
%Straight Lines [id:da6594475352061899] 
\draw    (198.43,90.71) -- (245.43,83.29) ;
%Shape: Free Drawing [id:dp020898871037853928] 
\draw  [line width=3] [line join = round][line cap = round] (309.43,153.71) .. controls (309.43,153.71) and (309.43,153.71) .. (309.43,153.71) ;
%Shape: Free Drawing [id:dp4349336473189229] 
\draw  [line width=3] [line join = round][line cap = round] (224.43,160.71) .. controls (224.43,160.71) and (224.43,160.71) .. (224.43,160.71) ;
%Shape: Free Drawing [id:dp16450328641071899] 
\draw  [line width=3] [line join = round][line cap = round] (260.43,79.71) .. controls (260.43,79.71) and (260.43,79.71) .. (260.43,79.71) ;

% Text Node
\draw (207,82.4) node [anchor=north west][inner sep=0.75pt]  [font=\small]  {$1$};
% Text Node
\draw (270.71,251.58) node [anchor=north west][inner sep=0.75pt]  [font=\small]  {$2$};
% Text Node
\draw (319,87.4) node [anchor=north west][inner sep=0.75pt]  [font=\small]  {$3$};
% Text Node
\draw (156,58.4) node [anchor=north west][inner sep=0.75pt]  [font=\small]  {$4$};
% Text Node
\draw (261,215.4) node [anchor=north west][inner sep=0.75pt]  [font=\small]  {$5$};
% Text Node
\draw (361,50.4) node [anchor=north west][inner sep=0.75pt]  [font=\small]  {$6$};
% Text Node
\draw (330,75.4) node [anchor=north west][inner sep=0.75pt]  [font=\scriptsize]  {$\textcolor[rgb]{1,0,0}{0}$};
% Text Node
\draw (353,50.4) node [anchor=north west][inner sep=0.75pt]  [font=\scriptsize]  {$\textcolor[rgb]{1,0,0}{0}$};
% Text Node
\draw (267.71,273.83) node [anchor=north west][inner sep=0.75pt]  [font=\scriptsize]  {$\textcolor[rgb]{1,0,0}{\epsilon _{2}}$};
% Text Node
\draw (259.71,203.4) node [anchor=north west][inner sep=0.75pt]  [font=\scriptsize]  {$\textcolor[rgb]{1,0,0}{\epsilon }\textcolor[rgb]{1,0,0}{_{5}}$};
% Text Node
\draw (207.71,97.4) node [anchor=north west][inner sep=0.75pt]  [font=\scriptsize]  {$\textcolor[rgb]{1,0,0}{\epsilon }\textcolor[rgb]{1,0,0}{_{1}}$};
% Text Node
\draw (182.71,56.4) node [anchor=north west][inner sep=0.75pt]  [font=\scriptsize]  {$\textcolor[rgb]{1,0,0}{\epsilon }\textcolor[rgb]{1,0,0}{_{4}}$};
% Text Node
\draw (347.71,90.4) node [anchor=north west][inner sep=0.75pt]  [font=\scriptsize]  {$\textcolor[rgb]{1,0,0}{l_{1}}$};
% Text Node
\draw (276.71,174.4) node [anchor=north west][inner sep=0.75pt]  [font=\scriptsize]  {$\textcolor[rgb]{1,0,0}{l_{2}}$};
% Text Node
\draw (308.71,105.4) node [anchor=north west][inner sep=0.75pt]  [font=\scriptsize]  {$\textcolor[rgb]{1,0,0}{l}\textcolor[rgb]{1,0,0}{_{3}}$};
% Text Node
\draw (292.71,208.4) node [anchor=north west][inner sep=0.75pt]  [font=\scriptsize]  {$\textcolor[rgb]{1,0,0}{l}\textcolor[rgb]{1,0,0}{_{4}}$};
% Text Node
\draw (242.71,175.4) node [anchor=north west][inner sep=0.75pt]  [font=\scriptsize]  {$\textcolor[rgb]{1,0,0}{a_{3}}$};
% Text Node
\draw (217.71,116.4) node [anchor=north west][inner sep=0.75pt]  [font=\scriptsize]  {$\textcolor[rgb]{1,0,0}{a}\textcolor[rgb]{1,0,0}{_{2}}$};
% Text Node
\draw (193.71,129.4) node [anchor=north west][inner sep=0.75pt]  [font=\scriptsize]  {$\textcolor[rgb]{1,0,0}{a}\textcolor[rgb]{1,0,0}{_{4}}$};
% Text Node
\draw (228.71,208.4) node [anchor=north west][inner sep=0.75pt]  [font=\scriptsize]  {$\textcolor[rgb]{1,0,0}{a}\textcolor[rgb]{1,0,0}{_{1}}$};
% Text Node
\draw (230.71,87.4) node [anchor=north west][inner sep=0.75pt]  [font=\scriptsize]  {$\textcolor[rgb]{1,0.01,0.01}{b_{1}}$};
% Text Node
\draw (315.71,56.4) node [anchor=north west][inner sep=0.75pt]  [font=\scriptsize]  {$\textcolor[rgb]{1,0.01,0.01}{b_{2}}$};
% Text Node
\draw (286.71,87.4) node [anchor=north west][inner sep=0.75pt]  [font=\scriptsize]  {$\textcolor[rgb]{1,0.01,0.01}{b}\textcolor[rgb]{1,0.01,0.01}{_{3}}$};
% Text Node
\draw (206.71,59.4) node [anchor=north west][inner sep=0.75pt]  [font=\scriptsize]  {$\textcolor[rgb]{1,0.01,0.01}{b_{4}}$};

\end{tikzpicture}
\caption{Length Diagram Corresponding to Configurations (1), (4), (5)}
\end{figure}
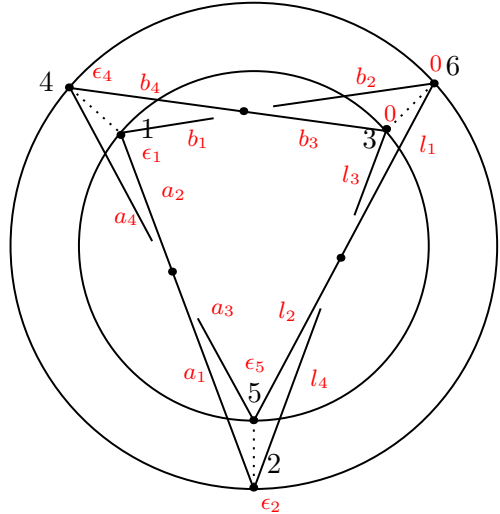
\end{center}
\par Now observe that if \(p_4\) is rotated along circle \(p_{2,4,6}\) more than the \(\epsilon\)
amount allotted by \(D'\), in either direction, then \(\frac{\tilde{l_1}}{\tilde{l_3}}\) will remain fixed. Assume for the sake of contradiction that \(\frac{\tilde{a_4}}{\tilde{a_2}}\) remains fixed and, without loss of generality, that the direction of rotation of \(p_4\) is counterclockwise. Since \(d\tilde{a_2}(p_4)\) is 0 along the \(p_4,p_5\) line and attains max/min values along the \(p_1,p_2\) parallel due to similar triangles and \(d\tilde{a_4}(p_4)\) is 0 along the parallel by similar triangles and attains max/min values along the \(p_4,p_5\) line, and the fact that \(d\tilde{a_2}(p_4), d\tilde{a_4}(p_4)\) vary trigonometrically with respect to angle, there is a single one-dimensional path spanned out within the envelope formed by the line between \(p_4\) and \(p_5\) and a line parallel to the \(p_1, p_2\) line centered at \(p_4\). Therefore, the tangent \(d\theta\) with respect to \(p_{2,4,6}\) 
can be decomposed into a component along the parallel line and along the \(p_4\),\(p_5\) line. Applying the law of sines to such components, the relationship \(\frac{\tilde{a_4}}{\tilde{a_2}}=\frac{\sin(\theta_1+\theta_2)}{\sin(\theta_1)} = c\) must be held for some constant \(c\) in order for the configuration to be \textit{bad}, where \(\theta_1\) is the angle between the oriented tangent to \(p_4\) and the ray \(\overrightarrow{p_5p_4}\) and \(\theta_2\) is the angle between \(\overrightarrow{p_2p_1}\) and \(\overrightarrow{p_5p_4}\). If we instead rotate \(p_4\) by \(\epsilon+\epsilon'<2\epsilon\), the LoS relationship fails to hold since \(\sin(\theta)\) is concave down when positive. Therefore, any \textit{bad} configuration with 
respect to the crossing rules may be taken to a \textit{good} configuration by the rotation of \(p_4 \in x \in \tilde{M}'\) by any \(\epsilon' < \epsilon\). Now, consider some \textit{good} configuration. Assume for sake of contradiction that there exists some \(\epsilon_0\) such that rotating \(p_4\) by \(\epsilon+\epsilon_0\) yields a \textit{bad} configuration. Then, by the above argument, rotating back to \(\epsilon\) will yield a \textit{good} configuration through all values within \([\epsilon-\epsilon_0,\epsilon+\epsilon_0]\). Therefore, there is always such an interval of \textit{good} configurations. $\boxed{}$

\hspace{4mm}
\par After having defined \textit{good} and \textit{bad} configurations for certain configuration types, we may now leverage Lemma 3.11 to prove that either certain configuration types must not exist in multiple planes, or there is a trefoil on \(\gamma\). If this can be done for every configuration type, then we will be finished as either there is a trefoil or no configuration type will exist in multiple planes, which contradicts the non-planarity of the entire \(S_I([D' \cap Q_\gamma])\) class given that the Conway quadratic term of \(\gamma\) is odd.
\begin{lemma}
    If \(\exists\) an interval \(I_\delta \in S_I([D' \cap Q_\gamma])\) of six strands of \(\gamma\) that consist of an \(\epsilon'\)-rotation about some line \(\ell\) of a type (1), (4), or (5) configuration \(x\), then \(\gamma\) has an inscribed trefoil.
\end{lemma}
\par \textit{Proof.} Assume that we have some \textit{bad} configuration \(x \in P\) for the \(\epsilon'\)-rotated configuration. Otherwise, the configuration \(x\) is \textit{good} and \(\gamma\) has an inscribed trefoil. Without loss of generality, assume that \(x\) is of type (1) and that \(\ell\) lies on \(p_{2,5}\). For all such \textit{bad} configurations, assume that \(\gamma \in P\) locally about \(p_2\) and \(p_5\). Otherwise, we may pick some orthogonal height along \(\gamma\) to bypass the \textit{bad} condition and find a trefoil. Then, we may find the tangent vector of the 1-dimensional path of \textit{bad} configurations through \(p_5\) when varying the point in the plane \(P\) with respect to all other points fixed. If we move \(p_5\) some infinitesimal amount \(\delta\) to \(p_5'\), then the vector \(\overrightarrow{p_5p_5'}\) may be componentized in terms of unit vectors \(\vec{v_0}, \vec{v_1}\), the unit vectors pointing along the same directions as \(\overrightarrow{p_1p_2}\) and \(\overrightarrow{p_4p_5}\), respectively. Letting, \(\overrightarrow{p_5p_5'} = \sum e_i\vec{v_i}\), a \textit{bad} configuration with \(p_5'\) instead of \(p_5\) must satisfy a relation: $$\frac{\frac{a_1-e_1\frac{a_4}{a_4+a_3+e_0}}{a_1+a_2}}{\frac{a_3+e_0}{a_3+a_4+e_0}}=\frac{\frac{a_1}{a_2+a_1}}{\frac{a_2}{a_3+a_4}}=\frac{\tilde{a_1}}{\tilde{a_3}}=c,$$ where \(c\) is some constant since none of the points \(p_1, p_3, p_4, p_6\) are 
perturbed. The relation is found by first moving along the \(\vec{v_0}\) direction \(e_0\), increasing \(a_3\) to \(a_3 + e_0\) and changing no other \(a_i\) lengths, then along the \(\vec{v_1}\) direction \(e_1\), which preserves \(\tilde{a_3}\) after the \(\vec{v_0}\) movement, decreases \(a_1\) by \(e_1 \frac{a_4}{a_4+a_3+e_0}\) by similar triangles, and preserves the sum length of \(a_1+a_2\). Simplifying this expression, the condition for the configuration with \(p_5'\) replacing \(p_5\) to be \textit{bad} becomes \(e_1 = -\frac{a_1}{a_3}e_0\). Therefore, \(\overrightarrow{p_5p_5'} = \sum e_i \vec{v_i} = e_0\vec{v_0} + e_1\vec{v_1} = \frac{e_0}{a_3}(a_3\vec{v_0}-a_1\vec{v_1})\) when \(p_5'\) is a point within a \textit{bad} configuration given that \(p_5\) is already within a \textit{bad} configuration. 
Considering these vectors with respect to the intersection point between \(p_1p_2\) and \(p_4p_5\), \(\overrightarrow{p_2p_5}||\overrightarrow{p_5p_5'}\), meaning that the 1-dimensional path lies exactly on the line \(p_{2,5}\).
\par Assuming \(\gamma \in P \cap p_{2,5}\) locally about \(p_5\),
we see that it is possible to follow \(\gamma\) in a unique direction from \(p_5\) such that the distance from the point \(p_2\) is minimized locally, as \(\gamma\) is smooth. Following this direction, there must exist a point \(p_5'\) that lies on a tubular neighborhood of the line segment from \(p_2\) to \(p_5\). Otherwise, \(\gamma\) will be non-injective or the segment \(\gamma_{5 \mapsto 6}\) taking \(p_5\) to \(p_6\) will intersect \(p_2\) first, a contradiction. Assume \(p_5' \in P\), as picking orthogonal heights for \(p_5' \notin P\) gives a trefoil. Then, \(p_5'\) lies outside of the 1-dimensional path of \(bad\) configurations \(x(t)=(p_1, p_2, p_3, p_4, p_5(t), p_6)\) with \(x(0)=x\) for \(p_5(t)=p_5\), which exists by Lemma 3.11, so \((p_1, p_2, p_3, p_4, p_5', p_6)\) is a \textit{good} configuration and there exists a trefoil on \(\gamma\). \(\boxed{}\)

\begin{center}
\begin{figure}
\tikzset{every picture/.style={line width=0.75pt}} %set default line width to 0.75pt        

\begin{tikzpicture}[x=0.75pt,y=0.75pt,yscale=-1,xscale=1]

%Straight Lines [id:da36651073474364315] 
\draw    (329,78) -- (312,214) ;
%Shape: Free Drawing [id:dp11136754995618992] 
\draw  [line width=3] [line join = round][line cap = round] (312.5,213) .. controls (312.5,213) and (312.5,213) .. (312.5,213) ;
%Shape: Free Drawing [id:dp9629881397109736] 
\draw  [line width=3] [line join = round][line cap = round] (328.5,78) .. controls (328.5,78) and (328.5,78) .. (328.5,78) ;
%Curve Lines [id:da40605120882010226] 
\draw    (312,214) .. controls (322,150) and (314,126) .. (354,96) ;
%Shape: Free Drawing [id:dp6632627892736123] 
\draw  [line width=3] [line join = round][line cap = round] (353.5,96) .. controls (353.5,96) and (353.5,96) .. (353.5,96) ;

% Text Node
\draw (318,201.4) node [anchor=north west][inner sep=0.75pt]    {$p_{5}$};
% Text Node
\draw (333,63.4) node [anchor=north west][inner sep=0.75pt]    {$p_{2}$};
% Text Node
\draw (356,91.4) node [anchor=north west][inner sep=0.75pt]    {$p_{5} '$};
% Text Node
\draw (332,121.4) node [anchor=north west][inner sep=0.75pt]    {$\gamma $};

\end{tikzpicture}
\caption{An example \(p_5'\) given \(p_2, p_5\).}
\end{figure}
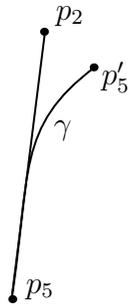
\end{center}

\begin{figure}
\begin{center}

\tikzset{every picture/.style={line width=0.75pt}} %set default line width to 0.75pt        

\begin{tikzpicture}[x=0.75pt,y=0.75pt,yscale=-1,xscale=1]

%Shape: Circle [id:dp5363648895882633] 
\draw   (36,134.1) .. controls (36,68.32) and (89.32,15) .. (155.1,15) .. controls (220.88,15) and (274.2,68.32) .. (274.2,134.1) .. controls (274.2,199.88) and (220.88,253.2) .. (155.1,253.2) .. controls (89.32,253.2) and (36,199.88) .. (36,134.1) -- cycle ;
%Shape: Circle [id:dp7508354147402772] 
\draw   (297,118) .. controls (297,75.47) and (331.47,41) .. (374,41) .. controls (416.53,41) and (451,75.47) .. (451,118) .. controls (451,160.53) and (416.53,195) .. (374,195) .. controls (331.47,195) and (297,160.53) .. (297,118) -- cycle ;
%Straight Lines [id:da019271384055250307] 
\draw    (91,33.33) -- (587,75.33) ;
%Straight Lines [id:da1020649498492634] 
\draw    (275,123.33) -- (587,75.33) ;
%Straight Lines [id:da34503358773486736] 
\draw    (168,253.33) -- (333,53.33) ;
%Straight Lines [id:da8724880167648956] 
\draw    (168,253.33) -- (587,75.33) ;
%Shape: Free Drawing [id:dp2865133157503654] 
\draw  [line width=3] [line join = round][line cap = round] (91.14,33.57) .. controls (91.14,33.57) and (91.14,33.57) .. (91.14,33.57) ;
%Shape: Free Drawing [id:dp5927115407067349] 
\draw  [line width=3] [line join = round][line cap = round] (169.14,251.57) .. controls (169.14,251.57) and (169.14,251.57) .. (169.14,251.57) ;
%Shape: Free Drawing [id:dp22208600503265807] 
\draw  [line width=3] [line join = round][line cap = round] (274.14,123.57) .. controls (274.14,123.57) and (274.14,123.57) .. (274.14,123.57) ;
%Shape: Free Drawing [id:dp02383483636352568] 
\draw  [line width=3] [line join = round][line cap = round] (332.14,54.57) .. controls (332.14,54.57) and (332.14,54.57) .. (332.14,54.57) ;
%Shape: Free Drawing [id:dp02784935039314318] 
\draw  [line width=3] [line join = round][line cap = round] (296.14,119.57) .. controls (296.14,119.57) and (296.14,119.57) .. (296.14,119.57) ;
%Shape: Free Drawing [id:dp3804972373775486] 
\draw  [line width=3] [line join = round][line cap = round] (449.14,133.57) .. controls (449.14,133.57) and (449.14,133.57) .. (449.14,133.57) ;
%Shape: Free Drawing [id:dp6153590629203671] 
\draw  [line width=3] [line join = round][line cap = round] (586.14,75.57) .. controls (586.14,75.57) and (586.14,75.57) .. (586.14,75.57) ;
%Straight Lines [id:da08727837390345328] 
\draw [color={rgb, 255:red, 255; green, 0; blue, 0 }  ,draw opacity=1 ]   (160.2,253.6) -- (333,53.33) ;
%Straight Lines [id:da7355444682272874] 
\draw [color={rgb, 255:red, 255; green, 0; blue, 0 }  ,draw opacity=1 ]   (286.25,109.25) -- (333,53.33) ;
%Straight Lines [id:da3186395597656715] 
\draw [color={rgb, 255:red, 255; green, 0; blue, 0 }  ,draw opacity=1 ]   (275,123.33) -- (449.33,133.67) ;
%Straight Lines [id:da7033046590975571] 
\draw [color={rgb, 255:red, 255; green, 0; blue, 0 }  ,draw opacity=1 ]   (160.2,253.6) -- (288.33,126.67) ;
%Straight Lines [id:da694633666762752] 
\draw [color={rgb, 255:red, 255; green, 0; blue, 0 }  ,draw opacity=1 ]   (292.25,122.25) -- (297,118) ;
%Straight Lines [id:da5184724443899007] 
\draw [color={rgb, 255:red, 255; green, 0; blue, 0 }  ,draw opacity=1 ]   (284.2,115.6) -- (297,118) ;
%Straight Lines [id:da4961418786055507] 
\draw [color={rgb, 255:red, 255; green, 0; blue, 0 }  ,draw opacity=1 ]   (275,123.33) -- (279.2,118.6) ;
%Straight Lines [id:da04903350440397758] 
\draw [color={rgb, 255:red, 255; green, 0; blue, 0 }  ,draw opacity=1 ]   (91,33.33) -- (279.2,111.6) ;
%Straight Lines [id:da35209239778419477] 
\draw [color={rgb, 255:red, 255; green, 0; blue, 0 }  ,draw opacity=1 ]   (91,33.33) -- (298.2,87.6) ;
%Straight Lines [id:da7677304555475017] 
\draw [color={rgb, 255:red, 255; green, 0; blue, 0 }  ,draw opacity=1 ]   (307,91.33) -- (449.33,133.67) ;
%Shape: Free Drawing [id:dp2467009673335736] 
\draw  [line width=3] [line join = round][line cap = round] (160.14,253) .. controls (160.14,253) and (160.14,253) .. (160.14,253) ;

% Text Node
\draw (90.5,36.73) node [anchor=north west][inner sep=0.75pt]    {$p_{1}$};
% Text Node
\draw (317,28.4) node [anchor=north west][inner sep=0.75pt]    {$p_{4}$};
% Text Node
\draw (170,256.73) node [anchor=north west][inner sep=0.75pt]    {$p_{5}$};
% Text Node
\draw (246,109.4) node [anchor=north west][inner sep=0.75pt]    {$p_{3}$};
% Text Node
\draw (299.67,95.9) node [anchor=north west][inner sep=0.75pt]    {$p_{6}$};
% Text Node
\draw (453,132.4) node [anchor=north west][inner sep=0.75pt]    {$p_{2}$};
% Text Node
\draw (593,68.4) node [anchor=north west][inner sep=0.75pt]    {$p$};
% Text Node
\draw (149,232.4) node [anchor=north west][inner sep=0.75pt]    {$\textcolor[rgb]{1,0,0}{0}$};
% Text Node
\draw (335,56.73) node [anchor=north west][inner sep=0.75pt]    {$\textcolor[rgb]{1,0,0}{0}$};
% Text Node
\draw (231,106.4) node [anchor=north west][inner sep=0.75pt]    {$\textcolor[rgb]{1,0,0}{0}$};
% Text Node
\draw (302,127.4) node [anchor=north west][inner sep=0.75pt]    {$\textcolor[rgb]{1,0,0}{+\epsilon }$};
% Text Node
\draw (452,149.4) node [anchor=north west][inner sep=0.75pt]    {$\textcolor[rgb]{1,0,0}{+\epsilon '\gg \epsilon }$};
% Text Node
\draw (81,56.4) node [anchor=north west][inner sep=0.75pt]    {$\textcolor[rgb]{1,0,0}{-\alpha \epsilon }$};
% Text Node
\draw (136,253.4) node [anchor=north west][inner sep=0.75pt]    {$p_{5} '$};

\end{tikzpicture}

\end{center}
\caption{A trefoil construction for configuration (3).}
\end{figure}
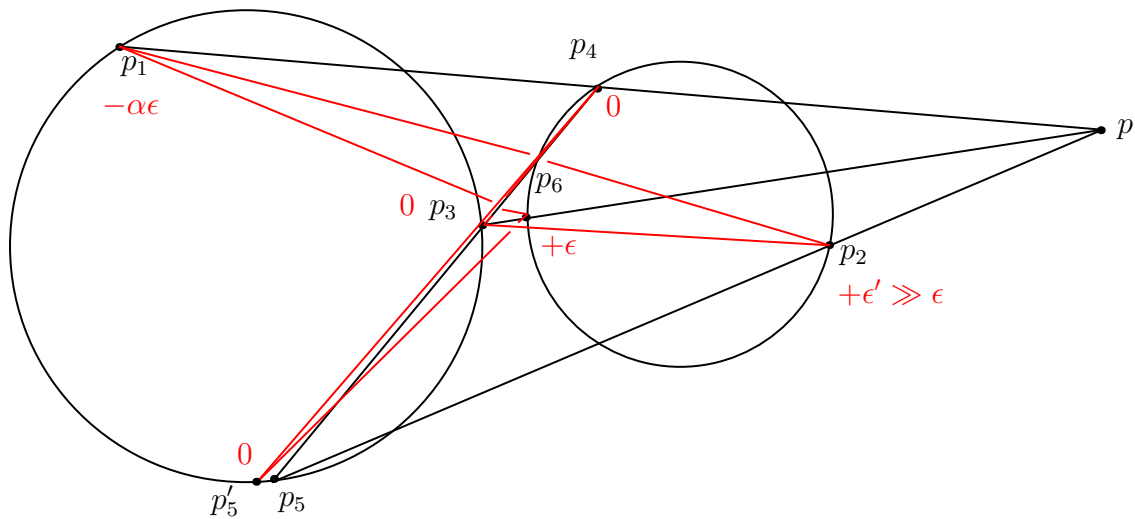

\begin{lemma}
    If \(\exists\) an interval \(I_\delta \in S_I([D' \cap Q_\gamma])\) of six strands of \(\gamma\) that consist of an \(\epsilon'\)-rotation about some line \(\ell\) of a type (3) configuration \(x\), then \(\gamma\) has an inscribed trefoil.
\end{lemma}
\par \textit{Proof.} Suppose \(x\) has the maximal number of points on \(\ell\): 3. If there were four points, then one of the circles \(p_{1,3,5}, p_{2,4,6}, p_{1,4,3,6}, p_{2,3,5,6}\) would be a line instead, which is impossible. Without loss of generality, let \(p_3, p_4, p_5\) be the three points. Note that all other cases reduce to this case, since having one or two points on \(\ell\) allows for at least one additional choice of orthogonal height along \(\gamma\).

Since \(p_6\) does not lie on the line \(p_3p_4p_5\), there exists some point \(p_5'\) within a tubular neighborhood of \(p_5\) not on the line, as the oriented strand \(\gamma_{5\mapsto6}\) must take \(p_5\) off of the line in order to intersect \(p_6\) before \(p_3\). Therefore, we may view the configuration in Figure 8 from the perspective of the plane passing through \(p_3, p_4,\) and \(p_5'\) while assigning orthogonal values to the other points with respect to rotation about the line \(p_{3,4,5}\). Picking \(\epsilon\) values for \(p_6, p_2\) as shown resolves the bottom-most crossing. Next, with the correct \(\alpha\) factor, the intersection between \(P\) and the line \(p_{1,6}\) may be chosen to lie between the lines \(p_{4,5}, p_{4,3}\) and give the desired crossing of the under-strand of the trefoil. The final crossings that must be resolved involve the \(p_{1,2}\) line, but so long as \(\epsilon'\) is chosen sufficiently, both crossings on this line may be made to pass over the lines \(p_{4,5}, p_{4,3}\) (shown is the case passing under these lines; note that either stick-knot is a trefoil). Thus, there is a trefoil on \(\gamma\). \(\boxed{}\)
\begin{lemma}
    If \(\exists\) an interval \(I_\delta \in S_I([D' \cap Q_\gamma])\) of six strands of \(\gamma\) that consist of an \(\epsilon'\)-rotation about some line \(\ell\) of a type (2) configuration \(x\), then \(\gamma\) has an inscribed trefoil.
\end{lemma}
\par \textit{Proof.} The last configuration we must construct a trefoil within is the (2) configuration with three of \(x\)'s points always colinear. Without loss of generality, consider a configuration in which \(p_2,p_4,p_6\) are colinear. Then, \(\exists p_2' \in \gamma\) on some tubular neighborhood of the ray \(\overrightarrow{p_4p_2}\), since \(\gamma_{2 \mapsto 3}\) must intersect \(p_3\) before \(p_4\) and therefore deviates from the line \(p_{2,4}\). Given that \(x \in I_\delta\), we may orient \(p_2'\) below an envelope of planes \(E_1\) containing a sub-interval \(I_{\delta'} \subset I_{\delta}\) from which orthogonal heights for \(p_i\) may be picked with respect to the plane \(p_{4,5,6}\).

\begin{figure}
\begin{center}
\tikzset{every picture/.style={line width=0.75pt}} %set default line width to 0.75pt        

\begin{tikzpicture}[x=0.75pt,y=0.75pt,yscale=-1,xscale=1]

%Straight Lines [id:da8680344831863158] 
\draw  [dash pattern={on 0.84pt off 2.51pt}]  (147.2,249.6) -- (490.2,249.6) ;
%Shape: Circle [id:dp9146012596056066] 
\draw  [dash pattern={on 0.84pt off 2.51pt}] (226,133.6) .. controls (226,81.91) and (267.91,40) .. (319.6,40) .. controls (371.29,40) and (413.2,81.91) .. (413.2,133.6) .. controls (413.2,185.29) and (371.29,227.2) .. (319.6,227.2) .. controls (267.91,227.2) and (226,185.29) .. (226,133.6) -- cycle ;
%Shape: Free Drawing [id:dp03202561289230155] 
\draw  [line width=3] [line join = round][line cap = round] (446.2,249.4) .. controls (446.2,249.4) and (446.2,249.4) .. (446.2,249.4) ;
%Shape: Free Drawing [id:dp930049161108002] 
\draw  [line width=3] [line join = round][line cap = round] (311.2,249.4) .. controls (311.2,249.4) and (311.2,249.4) .. (311.2,249.4) ;
%Shape: Free Drawing [id:dp5469197362604017] 
\draw  [line width=3] [line join = round][line cap = round] (184.2,249.4) .. controls (184.2,249.4) and (184.2,249.4) .. (184.2,249.4) ;
%Straight Lines [id:da15851972211946608] 
\draw  [dash pattern={on 0.84pt off 2.51pt}]  (365.2,22.4) -- (184,249) ;
%Straight Lines [id:da746120990712879] 
\draw  [dash pattern={on 0.84pt off 2.51pt}]  (365.2,22.4) -- (311.2,249.4) ;
%Straight Lines [id:da6797076948614693] 
\draw  [dash pattern={on 0.84pt off 2.51pt}]  (365.2,22.4) -- (446.2,249.4) ;
%Shape: Free Drawing [id:dp6599058346534736] 
\draw  [line width=3] [line join = round][line cap = round] (365.2,23.4) .. controls (365.2,23.4) and (365.2,23.4) .. (365.2,23.4) ;
%Shape: Free Drawing [id:dp29809589284156135] 
\draw  [line width=3] [line join = round][line cap = round] (239.2,180.4) .. controls (239.2,180.4) and (239.2,180.4) .. (239.2,180.4) ;
%Shape: Free Drawing [id:dp1309996466441834] 
\draw  [line width=3] [line join = round][line cap = round] (316.2,227.4) .. controls (316.2,227.4) and (316.2,227.4) .. (316.2,227.4) ;
%Shape: Free Drawing [id:dp29321854342686193] 
\draw  [line width=3] [line join = round][line cap = round] (411.2,151.4) .. controls (411.2,151.4) and (411.2,151.4) .. (411.2,151.4) ;
%Straight Lines [id:da057140746643310036] 
\draw    (296.6,230.6) -- (314.2,227.6) ;
%Straight Lines [id:da39873553212963353] 
\draw    (184,249) -- (268.2,214.6) ;
%Straight Lines [id:da6673282944556107] 
\draw    (411.2,150.6) -- (332.2,227.6) ;
%Straight Lines [id:da5189712180213879] 
\draw    (311.2,249.4) -- (238.2,180.4) ;
%Straight Lines [id:da3483533940759367] 
\draw    (446.2,249.4) -- (351.2,217.6) ;
%Straight Lines [id:da6045713058222302] 
\draw    (446.2,249.4) -- (314.2,227.6) ;
%Straight Lines [id:da756129642501689] 
\draw    (338.2,213.6) -- (238.2,180.4) ;
%Straight Lines [id:da7681814262169471] 
\draw    (327.2,233.6) -- (311.2,249.4) ;
%Straight Lines [id:da1439845841684786] 
\draw    (283.6,232.6) -- (184,249) ;
%Straight Lines [id:da7395388235752185] 
\draw    (277.2,209.6) -- (293.2,202.6) ;
%Straight Lines [id:da22267495845035534] 
\draw    (306.2,196.6) -- (411.2,150.6) ;

% Text Node
\draw (163,226.4) node [anchor=north west][inner sep=0.75pt]    {$p_{2}$};
% Text Node
\draw (293,250.4) node [anchor=north west][inner sep=0.75pt]    {$p_{4}$};
% Text Node
\draw (449,225.4) node [anchor=north west][inner sep=0.75pt]    {$p_{6}$};
% Text Node
\draw (299.6,203.1) node [anchor=north west][inner sep=0.75pt]    {$p_{1}$};
% Text Node
\draw (422,134.4) node [anchor=north west][inner sep=0.75pt]    {$p_{3}$};
% Text Node
\draw (215,166.4) node [anchor=north west][inner sep=0.75pt]    {$p_{5}$};
% Text Node
\draw (373,14.4) node [anchor=north west][inner sep=0.75pt]    {$p$};
% Text Node
\draw (313.2,252.8) node [anchor=north west][inner sep=0.75pt]  [font=\scriptsize]  {$\textcolor[rgb]{1,0,0}{{\textstyle 0}}$};
% Text Node
\draw (448.2,252.8) node [anchor=north west][inner sep=0.75pt]  [font=\scriptsize]  {$\textcolor[rgb]{1,0,0}{{\textstyle 0}}$};
% Text Node
\draw (420.2,155.8) node [anchor=north west][inner sep=0.75pt]  [font=\scriptsize]  {$\textcolor[rgb]{1,0,0}{+\epsilon _{3}}$};
% Text Node
\draw (236.2,165.4) node [anchor=north west][inner sep=0.75pt]  [font=\scriptsize]  {$\textcolor[rgb]{1,0,0}{{\textstyle 0}}$};
% Text Node
\draw (329.2,237) node [anchor=north west][inner sep=0.75pt]  [font=\scriptsize]  {$\textcolor[rgb]{1,0,0}{{\displaystyle +\epsilon _{1} \gg \epsilon _{3}}}$};
% Text Node
\draw (180,254.4) node [anchor=north west][inner sep=0.75pt]  [font=\scriptsize]  {$\textcolor[rgb]{1,0,0}{-\epsilon _{2} \gg \epsilon _{1}}$};

\end{tikzpicture}

\end{center}
\caption{A trefoil construction for configuration (2) with \(p_2'\) represented by \(p_2\).}
\end{figure}

\par Choosing orthogonal values as in Figure 9 (with \(p_2'\) represented by \(p_2\) for an arbitrarily-small tubular neighborhood) allows us to form a trefoil on \((p_1, p_2', p_3, p_4, p_5, p_6)\), in which the under-strand following the direction from \(p_5\) to \(p_6\) reaches over the over-strand from \(p_3\) to \(p_4\) as \(\epsilon_1\) is positive and far larger than \(\epsilon_3\). It reaches through the face formed by \(p_4,p_5,p_6\) and back around to \(p_3\), since \(\epsilon_2\) is negative and sufficiently larger in magnitude than either \(\epsilon_1\) or \(\epsilon_2\). Thus, given such a \(I_\delta\), \(\gamma\) has a trefoil as desired. \(\boxed{}\)
\newline
\par The last lemma we need is one that guarantees the existence of some \(I_\delta \in S_I([D' \cap Q_\gamma])\) for \(x \in I_\delta\) of any configuration type, given that \(S_I([D' \cap Q_\gamma])\) lies in at least two planes. Combining such a lemma with previous lemmas would guarantee a trefoil on \(\gamma\) given conditions on the Conway polynomial's quadratic term and \(\gamma\) smooth.

\begin{definition}
For a configuration type (i), let \(S_i\) be the set of configurations of that type.
\end{definition}

\begin{definition}
Let \(\theta(x)\) be the angle that configuration \(x\) makes with respect to the plane \(z=0\) in canonical coordinates. Then, we may alternatively define an \(\epsilon\)-rotation of a configuration \(x\) in \(S_I([D' \cap Q_\gamma])\) to be \(\{S_I(B_{\delta}(S_I^{-1}(x)) \cap [D' \cap Q_\gamma]) | \exists \delta(\epsilon) \hspace{1mm} \text{such that} \hspace{1mm} \forall x_0 \in S_I(B_{\delta}(S_I^{-1}(x)) \cap [D' \cap Q_\gamma]), |\theta(x_0)-\theta(x)|<\epsilon\}\).
\end{definition}

\begin{lemma}
If \(S_I([D' \cap Q_\gamma])\) lies in at least two planes, then \(\exists\) an interval \(I_\delta \in S_I([D' \cap Q_\gamma])\) of six strands of \(\gamma\) 
that consist of an \(\epsilon'\)-rotation about some line \(\ell\) of some configuration \(x\).
\end{lemma}
\par \textit{Proof.} Observe that for any \(x \in S_I([D' \cap Q_\gamma])\) of configuration types (1), (3), (4), or (5), there is an interval of \(S_I^{-1}(x)\) around which all points are of the same configuration type when taken back through \(S_I\). 
Therefore, \(S_{1,3,4,5}\) are open sets in \(S_I([D' \cap Q_\gamma])\). Taking the complement of their union, \(S_2\) is a closed set. Furthermore, \(\theta\) is a smooth function such that \(\epsilon\)-rotations are disjoint open intervals of \(S_I([D' \cap Q_\gamma])\) corresponding to elements with \(\theta' \neq 0\). Let such intervals be enumerated \(I_i\). Since there are at least two planes of \(S_I([D' \cap
Q_\gamma])\), there is at least \(I_1\) nonempty. Assume for sake of contradiction that \(\forall i, I_i \subset S_2\).
Otherwise, \(\exists x \in \cup_i I_i\) with \(x \in S_j\) for \(j\neq 2\), yielding an interval of an \(\epsilon\)-rotation. Then, \(\exists x \in \cup_i I_{i} \subsetneq S_2\), and there is a configuration interval of an \(\epsilon\)-rotation about \(x\) since it lies in the interior of \(S_2\). \(\boxed{}\)
\newline 
\newline
We now have a proof of Theorem 1.3. 

\textit{Proof.} Combine Lemmas 3.6, 3.17, and (3.12, 3.13, 3.14). \(\boxed{}\)

\section{Future Research}
\par It is of interest whether the property of having an inscribed trefoil is a knot invariant. Additionally, it is unknown whether other knot types may be inscribed in a class of smooth knots. We suspect that it is possible to find a figure-8 knot inscribed in some smooth class by considering a configuration manifold \(M_1\) with the data: 1) 3 parallel planes at various \(z\)-values on \(\partial B^4\) and a line \(\ell_i\) on each plane such that the projection down to one of the planes of each line is within an \(\epsilon\)-envelope of rotation and 2) some points \(\partial B^4 \cap \ell_i\) and a third point \(p_7 \in \partial B^4\) such that if \(\pi_z\) is the \(z\)-coordinate, \(\pi_z(p)_{p \in \ell_1} = \pi_z(p_7) + \epsilon_1\) for some sufficiently-small \(\epsilon_1\). Then, we would define \(M_1'\) as a perturbation of \(M_1\) such that each line has some specified \(\epsilon_i\)-rotation. However, we will then be working with a manifold of \(7\) points such that the algebra in Proposition 9 of \cite{hugelmeyer} must be modified. It is also of interest whether inscribed knots can always be reduced to trefoils by removing points. More specifically, we have the following conjecture.
\begin{conjecture}
If \(K\) is a smooth knot with odd quadratic term of its Conway polynomial, then every inscribed knot can be reduced to a trefoil by deleting some sequence of points and connecting the remaining segments in order. 
\end{conjecture}
\par Lastly, there may exist an alternative proof of Theorem 1.3 that involves proving the manifold of configurations of a certain type are transverse to \([D' \cap Q_\gamma]\).

\begin{theorem}
(0-Multijet Transversality Theorem) \cite{cantarella} Let \(M\) and \(N\) be smooth manifolds and let \(Z\) be a submanifold of \(C_n(J^0(M,N))\) for \(J^0(M,N) = M \times N\) the space of 0-jets. Given a smooth \(f:M \rightarrow N\), let the 0-jet of \(f\) be \(j^0f(\bold{p}) = (\bold{p}, f(\bold{p}))\) and the n-fold 0-multijet be \(j_n^0 f(\vec{p}) = (j^0f(p_1), \ldots, j^0f(p_n))\). Then, \(T_Z = \{f \in C^\infty(M,N) \hspace{1mm} | \hspace{1mm} j_n^0 f \pitchfork Z\}\) is \(C^m\)-dense in \(C^\infty(M,N)\) for any \(m\). In fact, if Z is compact, then \(T_Z\) is \(C^\infty\) open in \(C^\infty(M,N)\).
\end{theorem}
\begin{conjecture}
Let \(T\) be the submanifold of \(M'\) consisting of configurations in which one of the \((1,3,5), (4,6,2)\) triangles is a line under \(S_I\) projection. Then, \(T \pitchfork Q_\gamma\).
\end{conjecture}

We suspect that it may be possible to prove Conjecture 3.9 with the multijet transversality theorem. In this case, the \(T\)-configurations would be finite with respect to \(M'\) such that they would also be finite with respect to \(D\) or \(D'\). As a result, \([D' \cap Q_\gamma]\) would have finitely-many configurations consisting of three points on a line passing through \(S_I\), such that if the \(\epsilon\)-perturbation of the common plane \(P\) occurs at \((3)\), we would be able to pick some \(\epsilon' < \epsilon\) perturbation by the continuity of \(f\) such that the \(\epsilon'\)-perturbation could viewed as occurring to adjacent planar configurations \((1)\), \((3)\), or \((4)\) instead. Some adjustments have to be made to the multi-jet theorem before Conjecture 3.9 follows, since the configuration space \(C_6(S^3)\) has \(j^0_n f \cap Z\) a manifold of minimum dimension 18, whereas T is a submanifold of the 13-dimensional M.
\begin{figure}

\tikzset{every picture/.style={line width=0.75pt}} %set default line width to 0.75pt        

\begin{tikzpicture}[x=0.75pt,y=0.75pt,yscale=-1,xscale=1]

%Shape: Circle [id:dp39232467807252225] 
\draw   (241,150) .. controls (241,100.85) and (280.85,61) .. (330,61) .. controls (379.15,61) and (419,100.85) .. (419,150) .. controls (419,199.15) and (379.15,239) .. (330,239) .. controls (280.85,239) and (241,199.15) .. (241,150) -- cycle ;
%Straight Lines [id:da10479463020162427] 
\draw    (284,205) -- (345,177) ;
%Straight Lines [id:da3390912092287759] 
\draw    (332,189) -- (359,215) ;
%Straight Lines [id:da37212448222427974] 
\draw    (287,152) -- (322,181) ;
%Straight Lines [id:da3884674192877702] 
\draw    (287,152) -- (336,130) ;
%Straight Lines [id:da18231310997666772] 
\draw    (335,137) -- (340,141) ;
%Straight Lines [id:da38758601942114335] 
\draw    (290,100) -- (329,130) ;
%Straight Lines [id:da7331555993510881] 
\draw    (290,100) -- (377,113) ;
%Straight Lines [id:da10656035842842226] 
\draw    (336,111) -- (359,215) ;
%Straight Lines [id:da8107495871908967] 
\draw    (317,143) -- (307,162) ;
%Straight Lines [id:da5921675388472143] 
\draw    (348,147) -- (371,165) ;
%Straight Lines [id:da942435921283139] 
\draw    (356,172) -- (371,165) ;
%Straight Lines [id:da5314902766936644] 
\draw    (344,127) -- (377,113) ;
%Straight Lines [id:da1939522583564317] 
\draw    (302,173) -- (284,205) ;
%Straight Lines [id:da5091538082092653] 
\draw    (336,111) -- (327,123) ;
%Straight Lines [id:da7693063664453546] 
\draw    (323,129) -- (320,135) ;

% Text Node
\draw (282,37.4) node [anchor=north west][inner sep=0.75pt]    {$x\ \in M_{1} '\ \cap Q_{\gamma }$};

\end{tikzpicture}
\caption{Potential element \(x\) of \(M_1' \cap Q_\gamma\).}
\end{figure}
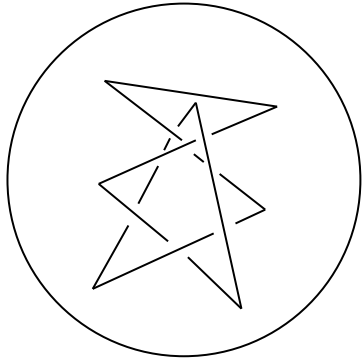
\section{Acknowledgements}
\par I would like to acknowledge Cole Hugelmeyer for his mentorship, ideas, and helpful discussions throughout this project. I would also like to thank the Caltech Summer Undergraduate Research Fellowship for funding and support, as well as Yi Ni for being my home institution mentor while I completed this SURF at Stanford University.
\bibliography{references}
\end{document}